\pdfoutput=1
\documentclass[11pt,a4paper]{article}

\usepackage[margin=1in]{geometry}
\usepackage[utf8]{inputenc}
\usepackage[T1]{fontenc}
\usepackage{lmodern}
\usepackage[numbers,sort&compress]{natbib}
\usepackage{hyperref}
\usepackage{url}
\usepackage{booktabs}
\usepackage{microtype}
\usepackage{xcolor}
\usepackage{amsmath}
\allowdisplaybreaks

\title{Near-Optimal Smooth Convex Optimization\\ with a Comparison Oracle}

\author{%
  Peter~Richt\'{a}rik\thanks{King Abdullah University of Science and Technology (KAUST), Thuwal, Saudi Arabia.
  Correspondence to: 
  \href{mailto:majiedammar.mahran@kaust.edu.sa}{\nolinkurl{majiedammar.mahran@kaust.edu.sa}}}
  \and
  Ammar~Mahran\footnotemark[1]%
}
\date{}

\usepackage{caption}
\usepackage{algorithm}
\usepackage{algpseudocode}
\usepackage{amsmath,amssymb,amsthm,mathtools}
\usepackage{aliascnt}
\usepackage{cleveref}
\usepackage{enumitem,etoolbox,needspace}
\usepackage{booktabs,tabularx,colortbl}

\theoremstyle{plain}
\newtheorem{theorem}{Theorem}[section]
\newaliascnt{proposition}{theorem}\newtheorem{proposition}[proposition]{Proposition}\aliascntresetthe{proposition}
\newaliascnt{lemma}{theorem}\newtheorem{lemma}[lemma]{Lemma}\aliascntresetthe{lemma}
\newaliascnt{corollary}{theorem}\newtheorem{corollary}[corollary]{Corollary}\aliascntresetthe{corollary}
\newtheorem{assumption}{Assumption}[section]
\theoremstyle{remark}
\newtheorem{remark}{Remark}[section]
\theoremstyle{plain}

\crefname{assumption}{Assumption}{Assumptions}
\crefname{lemma}{Lemma}{Lemmas}
\crefname{corollary}{Corollary}{Corollaries}
\crefname{proposition}{Proposition}{Propositions}
\crefname{remark}{Remark}{Remarks}

\definecolor{CompareTint}{HTML}{EDF5F7}

\newcommand{\norm}[1]{\lVert #1\rVert}
\newcommand{\R}{\mathbb{R}}
\newcommand{\ip}[2]{\langle#1,#2\rangle}
\DeclareMathOperator*{\argmin}{arg\,min}
\DeclareMathOperator{\spanof}{span}
\newcommand{\fa}{\qquad\text{for all }}
\newcommand{\Ball}{\mathbb B}
\newcommand{\Dir}{\mathsf{Direction}}
\newcommand{\SmallMin}{\mathsf{SmallMin}}

\BeforeBeginEnvironment{proof}{\Needspace{5\baselineskip}}
\newcommand{\pstep}[1]{\par\smallskip\noindent\emph{#1.}\enspace}
\newenvironment{proofsketch}{\begin{proof}[Proof sketch]}{\end{proof}}

\definecolor{linkcolor}{rgb}{0.45,0,0.05}
\hypersetup{
  pdftitle={Near-Optimal Smooth Convex Optimization with a Comparison Oracle},
  pdfauthor={Peter Richtarik, Ammar Mahran},
  colorlinks=true,
  linkcolor=linkcolor,
  citecolor=linkcolor,
  filecolor=linkcolor,
  urlcolor=linkcolor,
}

\begin{document}

\maketitle

\begin{abstract}
We study the minimization of a convex $d$-dimensional function with $L$-Lipschitz gradient when the only available information is a comparison oracle that reports which of two points has the smaller objective value.
Given an initial point within distance $R$ of a minimizer, the constants $L$ and $R$, and $A = LR^2/\varepsilon$, we give a deterministic method that returns an $\varepsilon$-optimal point using $O(\sqrt A\,[d\log(2dA)+\log^2(2A)])$ comparisons for $A\ge1$, counting all comparisons spent in inner subroutines.
The method is based on a normalized variant of Nemirovski's subspace method~\citep{nemirovski}:
each iteration performs a line search along an estimated gradient direction followed by a minimization over a two-dimensional affine subspace, and both subroutines rely on a comparison oracle and basic arithmetic operations only.
With exact function values, combining finite-difference gradients with an accelerated method gives the well-known bound of $O(d\sqrt A)$ evaluations~\citep{nesterovspokoiny2017,dorn2026}; 
our bound shows that one bit per query suffices to match it up to logarithmic factors.
Since two function values determine one comparison, the exact function-value lower bound of~\citet{dorn2026} implies that our bound is optimal up to logarithmic factors when $32\le A\le d^2$.

\end{abstract}



\section{Introduction}\label{sec:intro}

We consider the unconstrained convex optimization problem $\min_{x \in \R^d} f(x)$, where $f: \R^d \to \R$ denotes a $d$-dimensional objective function.
Throughout, we work under the following standing assumption.

\begin{assumption}\label{ass:setup}
The function $f:\R^d\to\R$ is convex, differentiable, and $L$-smooth with $L>0$, and it attains its minimum at some point $x^\star$ with $\norm{x^0-x^\star}\le R$ for a given initial point $x^0$ and some $R>0$:
\begin{equation}\label{eq:setup}
\begin{gathered}
 f(v)\ge f(u)+\ip{\nabla f(u)}{v-u} \fa u,v\in\R^d,\\
 \norm{\nabla f(v)-\nabla f(u)}\le L\norm{v-u} \fa u,v\in\R^d,\\
 x^\star\in\argmin_{u\in\R^d}f(u), \qquad \norm{x^0-x^\star}\le R .
\end{gathered}
\end{equation}
\end{assumption}

\paragraph{Target accuracy.}
A common goal in convex optimization is, for a given tolerance $\varepsilon > 0$, to find a point $x_{\rm out}$ with
\begin{equation}\label{eq:target}
 f(x_{\rm out})-f^\star \le \varepsilon .
\end{equation}
By smoothness (\Cref{ass:setup}), $f(x^0)-f^\star\le\frac L2\norm{x^0-x^\star}^2\le LR^2/2$, making $LR^2$ an a priori scale for the initial optimality gap.
We measure the target accuracy relative to this scale by $A\coloneqq LR^2/\varepsilon$, so that larger values of $A$ correspond to higher accuracy, and $1/A=\varepsilon/(LR^2)$ is the relative tolerance.
The quantity $A$ is dimensionless: replacing $f$ by $cf$ with $c>0$ multiplies $L$ and the gap $f(x)-f^\star$ by $c$, so $A$ is unchanged if $\varepsilon$ is multiplied by $c$ as well, and replacing $f$ by $f(\cdot/s)$ and $x^0$ by $sx^0$ with $s>0$ turns $L$ into $L/s^2$ and $R$ into $sR$, which leaves $LR^2$ and $A$ unchanged.
Henceforth, we state our results in terms of $A$, with the understanding that the tolerance in~\eqref{eq:target} is $\varepsilon = LR^2/A$.
We call a point $x_{\rm out}$ satisfying~\eqref{eq:target} $\varepsilon$-optimal.

\paragraph{Notation.}
\label{par:notation}
Throughout, $\ip\cdot\cdot$ and $\norm\cdot$ denote the Euclidean inner product and norm, and $\Ball(c,r)$ denotes the closed Euclidean ball with center $c$ and radius $r \ge 0$.
We write $f^\star \coloneqq f(x^\star)$ for the optimal objective value and $\Delta_k \coloneqq f(x^k)-f^\star$ for the optimality gap at iterate $x^k$.
Further, we define $g_k \coloneqq \nabla f(x^k)$ and $G_k \coloneqq \norm{g_k}$ as the gradient and gradient norm at iterate $x^k$.
The notation $\widetilde O(\cdot)$ hides factors that are polylogarithmic in the problem dimension $d$ and target accuracy $A$.

\paragraph{Oracle model.}
In \Cref{sec:finite}, we consider a method that accesses $f$ only through the comparison oracle
\begin{equation}\label{eq:oracle}
 O_f(u,v)=\mathbf1\{f(u)\ge f(v)\} ,
\end{equation}
which may be queried at any pair of points $u,v \in \R^d$ and returns a single bit.
The method may choose each query based on all previous answers, and it performs exact real arithmetic.
We track the number of these oracle calls in our analysis, including all calls made inside subroutines.

\paragraph{Acceleration is not immediate.}
With exact function values, the complexity of the class of functions satisfying \Cref{ass:setup} is known up to constants for $32\le A\le d^2$.
Finite differences allow approximation of the gradient up to a small error, and combining them with an accelerated gradient method gives an $\varepsilon$-optimal point after $O(d\sqrt A)$ function values, deterministically with coordinate differences~\citep[Corollary~2.3]{dorn2026} and in expectation with random directions~\citep[Section~6]{nesterovspokoiny2017}, which is optimal in this regime~\citep[Corollary~2.3]{dorn2026}.
Here, as for gradient methods~\citep{nesterov1983}, \emph{acceleration} refers to the improvement of the dependence on the accuracy from $A$ to $\sqrt A$.

A comparison oracle only reveals the sign of $f(u)-f(v)$, and two standard routes to acceleration are not immediately accessible.
First, the steps of accelerated gradient methods~\citep{nesterov1983} scale with the gradient norm, while comparisons cannot reveal it:
since $O_{cf}=O_f$ for every $c>0$, comparisons determine at most the direction of $\nabla f$, but not its magnitude.
Second, the standard reduction to strongly convex problems by regularization~\citep{allenzhu2016} is not available.
Accelerated methods for strongly convex functions, such as the comparison method of~\citet{lobanov}, would be applied to the regularized function $f_\mu(x)=f(x)+\frac\mu2\norm{x-x^0}^2$ with $\mu$ of order $\varepsilon/R^2$.
However, comparisons of $f$ do not determine comparisons of $f_\mu$:
the functions $f$ and $cf$ produce the same comparisons for every $c>0$, while $f_\mu$ and $cf+\frac\mu2\norm{\cdot-x^0}^2$ do not, in general.

The method we present in \Cref{sec:finite} circumvents these issues.
It only uses unit gradient directions, minimizations along rays, and minimizations over planes, and each of these can be carried out to sufficient accuracy via comparisons (\Cref{lem:direction,lem:small}).

\paragraph{Contributions.}
Without strong convexity, the best known comparison methods for the class of functions satisfying \Cref{ass:setup} use $\widetilde O(dA)$ comparisons (normalized descent) or $\widetilde O(d^2)$ comparisons (cutting planes)~\citep[Theorems~2 and~3]{zhangli2024}.
\citet[Section~1]{zhangli2024} pose as an open question whether comparison methods can be built on accelerated gradient methods, noting that comparisons reveal the direction of the gradient but not its norm.
We show that acceleration is indeed possible with comparisons alone, and that---ignoring logarithmic factors---it improves on both bounds when $1<A<d^2$:
\begin{enumerate}[leftmargin=*]
\item
\emph{An exact normalized subspace method} (\Cref{sec:ideal}).
We first analyze a normalized variant of Nemirovski's subspace method~\citep{nemirovski}, which only uses gradient directions, a line search, and a minimization over a two-dimensional affine subspace.
With exact subproblem solves, it guarantees $f(x^K)-f^\star \le 2LR^2 / K^2$ after $K$ iterations, so that $K=\lceil\sqrt{2A}\,\rceil$ iterations suffice for an $\varepsilon$-optimal point.
These iterations do not immediately translate to oracle calls in the comparison model, since each of them uses an exact gradient and two exact minimizations.
While the method is close to known subspace methods~\citep{narkiss,kuruzov}, its added value in the context of the present work lies in its proof, which serves as a template for a comparison-based method.

\item
\emph{An accelerated comparison method} (\Cref{sec:finite}).
We show that this guarantee survives errors in every quantity its proof relies on, and that every subroutine of the method can be implemented with the comparison oracle alone.
The resulting method returns an $\varepsilon$-optimal point after $O(\sqrt A[d\log(2dA)+\log^2(2A)]) = \widetilde O(d\sqrt A)$ comparisons, counting all comparisons made in its subroutines.
For $32 \le A\le d^2$, this is optimal up to logarithmic factors: 
since two function values determine one comparison, the exact function-value lower bound of~\citet{dorn2026} gives $\Omega(d\sqrt A)$ comparisons (\Cref{cor:lower-bound}).
\end{enumerate}
Combined with a comparison-based cutting-plane method for $A>d^2$ (\Cref{prop:cutting}), this determines the comparison complexity of the class up to logarithmic factors for $A\ge32$.
If each reply is correct with probability at least a known $p>1/2$ given the previous replies, majority votes preserve the accelerated rate with high probability (\Cref{cor:noisy}).
\paragraph{Related work.}
Methods that only compare function values have a long history in derivative-free optimization, among them Random Pursuit~\citep{stich2013}, which minimizes along random directions and is invariant under strictly increasing transformations of the objective, the stochastic three points method~\citep{bergou2020}, and Gradientless Descent~\citep{golovin2020}.
For smooth convex functions, these methods need of order $dLR^2/\varepsilon$ iterations up to logarithmic factors, with $R$ a radius of the level set $\{f\le f(x^0)\}$, or assume strong convexity; 
an accelerated variant of Random Pursuit is proposed in~\citet{stich2013} as a heuristic without guarantees.
\citet{jamieson2012} study optimization with noisy comparisons and with noisy function values for strongly convex functions, with rates that depend on how fast the probability of a correct reply approaches $1/2$ as the two function values approach each other, as in the Thurstone~\citep{thurstone1927} or Bradley--Terry~\citep{bradleyterry1952} models.
Their lower bounds~\citep[Theorem~1]{jamieson2012}, and that of~\citet[Theorem~14]{kumagai2017} for Gaussian noise, restrict the queries to a compact set, so they do not apply directly to the unrestricted queries in~\eqref{eq:oracle}.
Comparison-based optimization was further studied by~\citet{karabag} and~\citet{zhangli2024}.
For smooth convex functions and comparisons that are correct with a fixed probability $p>1/2$, a model that goes back to~\citet{feige1994}, \citet{saha} obtain $\widetilde O(dA/(p-1/2)^2)$ comparisons.
Comparison-based estimation of gradient directions appears in~\citep{cai,zhangli2024,tao}.
\citet{lobanov} obtain acceleration for strongly convex functions with a randomized coordinate method: 
their accelerated analysis treats the one-dimensional subproblems as solved exactly~\citep[Appendix~C.2]{lobanov}.
\citet{scheinberg} study comparison oracles without an underlying objective, and measure accuracy by the distance from a preference level set to the optimal solutions, under convexity, regularity, and growth conditions.
For exact function-value oracles, \citet{wu2026} and~\citet{dorn2026} give lower bounds that match the finite-difference upper bound $O(d\sqrt A)$ in a moderate-accuracy regime, for deterministic methods with queries in a ball and for randomized methods with unrestricted queries, respectively.
Concurrently with this work, the automated research framework LongTails lists a deterministic comparison method for the same class with $O((1+\sqrt A)(d\ell+\ell^3))$ comparisons, where $\ell=1+\log(d(1+A))$~\citep[Appendix~C, Table~32]{tyurin2026longtails}.
The method, only given in an associated repository~\citep{longtails}, estimates gradient directions as in~\citet[Algorithm~2]{zhangli2024} and, in each iteration, minimizes $f$ approximately over a three-dimensional box around $x^0$ spanned by $x^k-x^0$, a new normalized direction estimate, and the sum of the previous ones.
Our method uses one- and two-dimensional searches instead, and our bound has the additive term $\log^2(2A)$ in place of $\ell^3$; this changes the order of the bound only when $d$ is of smaller order than $\ell^2$, since otherwise the term of order $d\ell$, present in both bounds, dominates.

\Cref{tab:comparison} summarizes the complexity bounds discussed above.

\begin{table}[t]
\centering
\caption{Query complexity of minimizing a function that satisfies \Cref{ass:setup}, with $A=LR^2/\varepsilon$; $\widetilde O$ hides logarithmic factors.
For $32\le A\le d^2$, our method matches, up to logarithmic factors, both the comparison lower bound and the complexity with exact function values.}
\label{tab:comparison}
\small
\renewcommand{\arraystretch}{1.25}
\begin{tabularx}{\textwidth}{@{}>{\raggedright\arraybackslash}p{.34\textwidth}
 >{\centering\arraybackslash}p{.2\textwidth}
 >{\raggedright\arraybackslash}X@{}}
\toprule
\textbf{Method} & \textbf{Queries} & \textbf{Notes}\\
\midrule
\multicolumn{3}{@{}l}{\emph{Comparison oracle, upper bounds}}\\
Ellipsoid method\newline{\footnotesize\citet[Thm.~2]{karabag}}
& $\widetilde O(d^4)$ & Compact domain; Lipschitz $f$.\\
Normalized descent\newline{\footnotesize\citet[Thm.~2]{zhangli2024}}
& $\widetilde O(dA)$ & \\
Cutting planes\newline{\footnotesize\citet[Thm.~3]{zhangli2024}}
& $\widetilde O(d^2)$ & Self-contained proof with $O(d^2\log A\log(dA))$ in \Cref{prop:cutting}.\\
Accelerated coordinate search\newline{\footnotesize\citet[Thm.~4.1]{lobanov}}
& $\widetilde O(d\sqrt\kappa)$ & Requires $\mu$-strong convexity, $\kappa=L/\mu$; exact line searches; randomized.\\
\rowcolor{CompareTint}
Comparison-based normalized subspace method\newline{\footnotesize This paper, \Cref{thm:finite}}
& $\boldsymbol{\widetilde O(d\sqrt A)}$ & \textbf{Optimal up to logarithmic factors for $\boldsymbol{32\le A\le d^2}$}.\\
\midrule
\multicolumn{3}{@{}l}{\emph{Comparison oracle, lower bound}}\\
Any method\newline{\footnotesize \Cref{cor:lower-bound}, via \citet[Thm.~2.1]{dorn2026}}
& $\Omega(d\min\{d,\sqrt A\})$ & $A\ge32$; two function values determine one comparison.\\
\midrule
\multicolumn{3}{@{}l}{\emph{Exact function-value oracle, for reference}}\\
Finite differences and acceleration\newline{\footnotesize\citet[Cor.~2.3]{dorn2026}}
& $O(d\sqrt A)$ & Optimal for $32\le A\le d^2$.\\
\bottomrule
\end{tabularx}
\par\smallskip
\begin{minipage}{\textwidth}\footnotesize
The comparison upper bounds are deterministic unless noted otherwise, and they count all comparisons, including those made in inner subroutines.
The exception is the accelerated coordinate search: its analysis treats the line searches as exact~\citep[Appendix~C.2]{lobanov}, and it bounds the expected error; we state it for coordinate smoothness constants $L_i\le L$.
\end{minipage}
\end{table}

\section{Exact normalized subspace method}\label{sec:ideal}

We begin with an iterative method for minimizing $f$ that has access to normalized gradients and to exact minimizers over rays and planes.
Starting from the initial point $x^0$ of \Cref{ass:setup}, the method maintains a memory vector, initialized as $m^0 = 0 \in \R^d$, that accumulates unit gradients over time:
\begin{equation}\label{eq:ideal}
 m^{k+1}=m^k+n^k, \qquad n^k=g_k/G_k, \qquad\text{when }G_k>0 .
\end{equation}
If $G_k = 0$, then $x^k$ is already a minimizer of $f$, and the method stops.
Otherwise, it searches the steepest-descent ray from $x^k$,
\begin{equation}\label{eq:idealray}
 \alpha_k\in\argmin_{\alpha\ge0}f(x^k-\alpha n^k),\qquad
 p^{k+1}=x^k-\alpha_kn^k,
\end{equation}
so that $p^{k+1}$ is a minimizer of $f$ over the ray that starts at $x^k$ and points in the direction $-n^k$.
Next, the method minimizes $f$ over the affine subspace through the initial point $x^0$ that is spanned by the vector $p^{k+1} - x^0$ and the memory vector $m^{k+1}$:
\begin{equation}\label{eq:idealplane}
 V_{k+1}=\spanof\{p^{k+1}-x^0,m^{k+1}\},\qquad
 x^{k+1}\in\argmin_{x\in x^0+V_{k+1}}f(x).
\end{equation}
The subspace $V_{k+1}$ may have dimension less than two.
\Cref{alg:ideal} summarizes the method.

\begin{algorithm}[tbp]
\caption{Exact normalized subspace method}\label{alg:ideal}
\begin{algorithmic}[1]
\Require $x^0$, horizon $K\ge1$; exact gradients and attained exact ray and plane minima.
\State $m^0=0$.
\For{$k=0,\ldots,K-1$}
\State Compute $g_k=\nabla f(x^k)$ and $G_k=\norm{g_k}$.\label{line:ideal-gradient}
\If{$G_k=0$} \State \Return $x^k$. \EndIf
\State Set $n^k=g_k/G_k$ and $m^{k+1}=m^k+n^k$.
\State Solve the ray problem~\eqref{eq:idealray} for $p^{k+1}$.\label{line:ideal-ray}
\State Solve the plane problem~\eqref{eq:idealplane} for $x^{k+1}$.\label{line:ideal-plane}
\EndFor
\State \Return $x^K$.
\end{algorithmic}
\end{algorithm}

\begin{theorem}[Exact normalized subspace method]\label{thm:ideal}
Let \Cref{ass:setup} hold, let $K\ge1$, and suppose that the minima in~\eqref{eq:idealray} and~\eqref{eq:idealplane} are attained in every iteration of \Cref{alg:ideal}.
Then the point $x_{\rm out}$ returned by \Cref{alg:ideal} satisfies
\begin{equation}\label{eq:idealrate}
 f(x_{\rm out})-f^\star\le\frac{2LR^2}{K^2} .
\end{equation}
\end{theorem}

\Cref{alg:ideal} uses the gradient only through its direction $n^k$ and the test $G_k=0$, and it uses neither the value of the smoothness constant, $L$, nor that of the initial distance to a minimizer, $R$.
\begin{proof}[Proof of \Cref{thm:ideal}]
If \Cref{alg:ideal} stops at some $x^k$ with $G_k=0$, the returned point is a minimizer and~\eqref{eq:idealrate} holds.
We may therefore assume that $G_k>0$ for all $0 \le k < K$, so that $x_{\rm out}=x^K$.

\emph{Sufficient decrease.}
By $L$-smoothness of $f$, we have $f(x+s) \le f(x) + \ip{\nabla f(x)}{s} + \frac{L}{2} \norm{s}^2$ for all $x,s \in \R^d$.
Note that the gradient step $x^k - g_k / L$ lies on the ray in~\eqref{eq:idealray}, at $\alpha = G_k/L$, and $p^{k+1}$ lies in the affine subspace $x^0+V_{k+1}$ over which $x^{k+1}$ is optimal in~\eqref{eq:idealplane}.
Therefore,
\begin{equation}\label{eq:idealdescent}
  f(x^{k+1}) \le f(p^{k+1}) \le f(x^k-g_k/L) \le f(x^k) - \frac{G_k^2}{2L} ,
  \quad\text{i.e.,}\quad
  \Delta_k - \Delta_{k+1} \ge \frac{G_k^2}{2L} .
\end{equation}
In particular, the gaps $\Delta_k$ are non-increasing.
Now, if $\Delta_K=0$, the returned point is a minimizer and~\eqref{eq:idealrate} holds.
We may therefore assume that $\Delta_K>0$.
Then, $\Delta_k \ge \Delta_K > 0$ for all $0 \le k \le K$.

\emph{Orthogonality.}
For $k \ge 1$, the point $x^k$ minimizes the convex differentiable function $f$ over the affine subspace $x^0+V_k$, so the gradient $g_k$ is orthogonal to $V_k$.
Both $m^k$ and $x^k-x^0$ belong to $V_k$, and for $k=0$ both vectors are zero.
Therefore, for all $0\le k<K$,
\begin{equation}\label{eq:ideal-orth}
  \ip{g_k}{m^k} = 0, \qquad \ip{g_k}{x^k-x^0} = 0 .
\end{equation}

\emph{Norm of the memory vector.}
By~\eqref{eq:ideal} and~\eqref{eq:ideal-orth}, $\norm{m^{k+1}}^2 = \norm{m^k}^2 + 2 \ip{n^k}{m^k} + \norm{n^k}^2 = \norm{m^k}^2 + 1$, so $\norm{m^K}^2 = K$.

\emph{Upper bound on normalized gaps.}
Convexity and the second relation in~\eqref{eq:ideal-orth} give $\Delta_k \le \ip{g_k}{x^k-x^\star} = \ip{g_k}{x^0-x^\star}$.
Dividing by $G_k$, summing over $k$, and using~\eqref{eq:ideal}, the Cauchy--Schwarz inequality, and \Cref{ass:setup}, we obtain
\begin{equation}\label{eq:ideal-upper-sum}
 \sum_{k=0}^{K-1}\frac{\Delta_k}{G_k}
 \le\ip{m^K}{x^0-x^\star}
 \le\norm{m^K}\norm{x^0-x^\star}
 \le R\sqrt K .
\end{equation}

\emph{Lower bound on normalized gaps.}
Define $D_k \coloneqq \Delta_k-\Delta_{k+1}$ and $t_k \coloneqq D_k/\Delta_k^2$, both positive by~\eqref{eq:idealdescent}.
Then~\eqref{eq:idealdescent} reads $G_k\le\sqrt{2LD_k}$, so $\Delta_k/G_k\ge t_k^{-1/2}/\sqrt{2L}$.
Since $\Delta_{k+1} \le \Delta_k$, the sum of the $t_k$ is bounded by a telescoping sum:
\begin{equation}
  \sum_{k=0}^{K-1}t_k \le \sum_{k=0}^{K-1}\frac{D_k}{\Delta_k\Delta_{k+1}}
  = \sum_{k=0}^{K-1}\left(\frac1{\Delta_{k+1}}-\frac1{\Delta_k}\right)
 \le \frac1{\Delta_K}.
\end{equation}
Hence, by Jensen's inequality for the convex function $t\mapsto t^{-1/2}$,
\begin{equation}\label{eq:ideal-lower-sum}
 \sqrt{2L} \sum_{k=0}^{K-1}\frac{\Delta_k}{G_k} \ge\sum_{k=0}^{K-1} t_k^{-1/2} \ge K^{3/2} \Bigl(\sum_{k=0}^{K-1}t_k\Bigr)^{-1/2} \ge K^{3/2}\sqrt{\Delta_K}.
\end{equation}

\emph{Conclusion.}
Combining~\eqref{eq:ideal-upper-sum} and~\eqref{eq:ideal-lower-sum} gives $K^{3/2}\sqrt{\Delta_K}\le R\sqrt{2LK}$, that is, $\Delta_K \le 2LR^2/K^2$, which is \eqref{eq:idealrate}.
\end{proof}

\paragraph{Relation to existing methods.}
\Cref{alg:ideal} belongs to the family of subspace methods that go back to Nemirovski~\citep{nemirovski}.
The sequential subspace optimization (SESOP) method of~\citet{narkiss} minimizes $f$, in each iteration, over the affine subspace through $x^0$ spanned by the current gradient, the displacement $x^k-x^0$, and a weighted sum $\sum_iw_ig_i$ of the gradients at the previous iterates, with weights $w_0=1$ and $w_i=\frac12+(\frac14+w_{i-1}^2)^{1/2}$, so that $w_i\approx i/2$~\citep[Sections~2.1 and~2.3]{narkiss}.
Its analysis uses that the current gradient is orthogonal to $x^k-x^0$ and to the weighted sum of the previous gradients~\citep[Section~2.3]{narkiss}.
\Cref{alg:ideal} uses the same two orthogonality relations, but with the weights $w_i=1/G_i$ in place of $w_i\approx i/2$: 
its memory vector is now $m^k=\sum_{i=0}^{k-1} g_i/G_i$.
Its plane $x^0+V_{k+1}$ lies in the affine subspace through $x^0$ spanned by $x^k-x^0$, $g_k$, and $m^k$, which is the SESOP subspace with these weights, since $p^{k+1}-x^0=(x^k-x^0)-\alpha_kn^k$ and $m^{k+1}=m^k+n^k$; 
the current gradient enters through the line search rather than as a separate search direction.
\citet[Algorithm~2]{kuruzov}, a modification of Nemirovski's conjugate gradient method (proposed, with restarts, for strongly convex problems by~\citet[Section~7.3]{nemirovskyyudin1983}), is closer in structure to \Cref{alg:ideal}: it alternates a minimization over the plane through $x^0$ spanned by $x^k-x^0$ and the sum of all previous gradients with a gradient step with step size $1/(2L)$, and its memory vector sums the gradients rather than their directions.
The proof of \Cref{thm:ideal} follows the analysis of the ORTH method~\citep[Section~2.3]{narkiss}, which bounds a weighted sum of the gaps $\Delta_k$ through the two orthogonality relations; we use the weights $1/G_k$ instead (cf.\ the $1/\norm{\nabla f}$-weighted averaging of~\citet{levy2017}), for which $\norm{m^K}^2=K$, and the sufficient decrease enters through the lower bound~\eqref{eq:ideal-lower-sum}.
In contrast to these methods, neither step of \Cref{alg:ideal} uses the gradient norm, which a comparison oracle could not reveal since rescaling $f$ preserves every comparison but rescales $\nabla f$.
More generally, the iterates of \Cref{alg:ideal} do not change if $f$ is replaced by $\phi\circ f$ for a differentiable function $\phi$ with $\phi'>0$, since gradient directions and exact minimizers do not change; comparisons have the same invariance.

\paragraph{Attainment of the subproblem minima.}
\Cref{thm:ideal} assumes that the minima in~\eqref{eq:idealray} and~\eqref{eq:idealplane} are attained.
This need not be the case, even if $f$ has a minimizer.
The following example shows this for the ray problem.

\begin{remark}[A ray problem without a minimizer]\label{rem:attainment}
We consider the function
\begin{equation*}
   \begin{aligned}
  f(x) &= H(\psi(x)) , \\
  \psi(x) &= \varphi(x_1 - x_2) - x_2 , \\
  \varphi(s) &= \sqrt{1 + s^2} - s ,
 \end{aligned}
 \qquad\qquad
 H(s) = \begin{cases}
  0 , & s \le 0 , \\
  s^2/2 , & 0 < s < 1 , \\
  s - 1/2 , & s \ge 1 .
 \end{cases}
\end{equation*}
Here, $\varphi$ is convex, positive, and decreasing, with $|\varphi'| \le 2$ and $0<\varphi'' \le 1$.
Similarly, $H$ is convex and non-decreasing, with $1$-Lipschitz derivative taking values in $[0,1]$.
Since $\psi$ is convex and $H$ is convex and non-decreasing, $f$ is convex.
Moreover, $\nabla\psi(x)=(\varphi'(s),-\varphi'(s)-1)$ with $s=x_1-x_2$, so $\norm{\nabla\psi}\le\sqrt5$ and $\nabla\psi$ is $2$-Lipschitz.
Since $\nabla f=H'(\psi)\nabla\psi$, $|H'|\le1$, $H'$ is $1$-Lipschitz, and $\psi$ is $\sqrt5$-Lipschitz, the product rule gives
$\norm{\nabla f(x)-\nabla f(y)}\le1\cdot2\norm{x-y}+\sqrt5\cdot\sqrt5\norm{x-y}$, so $\nabla f$ is $7$-Lipschitz.
The minimum value $0$ is attained over an unbounded region including the point $x^\star = (3/2, 1/2)$.
The function $f$ therefore meets \Cref{ass:setup} with $x^0 = (-1,-1)$, $L=7$, and $R=3$ (the chosen values for $L,R$ are not the smallest possible).
Now, $\nabla f(x^0)=(-1,0)$ and $n^0=(-1,0)$, so $f(x^0-\alpha n^0)=f(-1+\alpha,-1)=\varphi(\alpha)+\frac{1}{2}$ for all $\alpha \ge 0$, which approaches $1/2$ as $\alpha \to \infty$, but never attains this limit.
Thus the ray problem~\eqref{eq:idealray} has no solution.
\end{remark}

A sufficient condition for the minima in~\eqref{eq:idealray} and~\eqref{eq:idealplane} to be attained is that the set of minimizers of $f$ is bounded.
For a convex function, this implies that all level sets are bounded~\citep[Corollary~8.7.1]{rockafellar1970}.
The ray and the plane then intersect the level sets $\{f\le f(x^k)\}$ and $\{f\le f(p^{k+1})\}$, respectively, in nonempty compact sets, over which the continuous function $f$ attains its minimum.

\section{Comparison-based normalized subspace method}\label{sec:finite}

\Cref{alg:ideal} relies on three exact computations: a normalized gradient, a minimizer over a ray, and a minimizer over a plane.
In this section, we present a method that mirrors \Cref{alg:ideal} but only has access to the comparison oracle~\eqref{eq:oracle}.

Note that none of the three computations can be carried out exactly with comparisons:
finitely many bits determine a direction or a minimizer only approximately, and the unrestricted minima need not be attained (\Cref{rem:attainment}).
We therefore
(i) estimate directions,
(ii) solve the ray and plane problems inexactly over large balls around $x^0$,
(iii) compare the output of each solver with the point it started from and keep the better one, which makes the objective values monotone, and
(iv) run for a number of iterations fixed in advance.
The last point replaces the stopping test $G_k=0$ of \Cref{alg:ideal}: the method has no test for a vanishing gradient and does not certify $\varepsilon$-optimality.

We proceed in three steps.
\Cref{sec:inexact} shows that the convergence proof of \Cref{thm:ideal} tolerates errors.
\Cref{sec:primitives} constructs comparison routines for the three computations, and \Cref{sec:algorithm} chooses their accuracies so that the errors stay within the tolerances of \Cref{sec:inexact}.

\subsection{Acceleration under relaxed conditions}\label{sec:inexact}

The proof of \Cref{thm:ideal} relies on four properties of the iterates:
the memory vector is a sum of unit gradients;
the gradient $g_k$ is orthogonal to the memory vector $m^k$;
it is orthogonal to the displacement $x^k-x^0$;
and each iteration decreases the objective by at least $G_k^2/(2L)$.
With comparisons, the unit gradients are only approximate, the two orthogonality relations hold only approximately, and the decrease is smaller.
The following theorem relaxes all four properties: the inequalities in~\eqref{eq:robust-conditions} bound the direction error, the products with the memory vector and with the displacement, and the decrease, and they only need to hold while the current gap exceeds the target $\varepsilon$.
We refer to them as the direction~\eqref{eq:cond-direction}, memory~\eqref{eq:cond-memory}, anchor~\eqref{eq:cond-anchor}, and descent~\eqref{eq:cond-descent} conditions.

\begin{theorem}[Acceleration under relaxed conditions]\label{thm:inexact-certificate}
Let $f:\R^d\to\R$ be convex and differentiable, with a minimizer $x^\star$.
Let $K\ge1$ be an integer, let $\varepsilon\ge0$, and let $x^0,\ldots,x^K\in\R^d$ be points whose optimality gaps are non-increasing, $\Delta_{k+1}\le\Delta_k$, and such that $\norm{x^0-x^\star}\le R$ for some $R>0$.
Let $u^0,\ldots,u^{K-1}$ be unit vectors, and define the memory vectors by $m^0=0$ and $m^{k+1}=m^k+u^k$ for $0\le k<K$.
Suppose that there are a constant $\eta>0$ and nonnegative numbers $\vartheta_k,r_k,\delta_k$ such that, for every $0\le k<K$ with $\Delta_k>\varepsilon$,
\begin{subequations}\label{eq:robust-conditions}
\begin{align}
 G_k>0\quad\text{and}\quad\norm{u^k-g_k/G_k}&\le\delta_k,\label{eq:cond-direction}\\
 \ip{u^k}{m^k}&\le\vartheta_k,\label{eq:cond-memory}\\
 \ip{g_k}{x^k-x^0}/G_k&\le r_k,\label{eq:cond-anchor}\\
 \Delta_k-\Delta_{k+1}&\ge\eta G_k^2.\label{eq:cond-descent}
\end{align}
\end{subequations}
With $S_K\coloneqq R\,(K+2\sum_{k=0}^{K-1}\vartheta_k)^{1/2}+\sum_{k=0}^{K-1}(r_k+R\delta_k)$, we then have
\begin{equation}\label{eq:robust-rate}
 \Delta_K\le\max\left\{\varepsilon,\frac{S_K^2}{\eta K^3}\right\}.
\end{equation}
\end{theorem}

The proof is given in \Cref{app:ideal}; it follows that of \Cref{thm:ideal}, and the errors accumulate in the quantity $S_K$.

For the iterates of \Cref{alg:ideal}, we have $u^k=n^k$ and $\vartheta_k=r_k=\delta_k=0$ by~\eqref{eq:ideal-orth}, and~\eqref{eq:idealdescent} holds with $\eta=1/(2L)$; with $\varepsilon=0$, this gives $S_K=R\sqrt K$, and~\eqref{eq:robust-rate} recovers~\eqref{eq:idealrate}.
More generally, the rate $O(LR^2/K^2)$ is retained above the target as long as $S_K=O(R\sqrt K)$ and $\eta$ is of order $1/L$.
This is the case if the direction errors $\delta_k$ and the normalized anchor products $r_k/R$ are of order $1/\sqrt K$ on average, and the memory products $\vartheta_k$ are, on average, of order one.

\subsection{Comparison primitives}\label{sec:primitives}

Recall that \Cref{alg:ideal} uses only gradient directions, not magnitudes.
Such directions can be estimated with finite differences, and the sign of a finite difference may be obtained via the comparison oracle \eqref{eq:oracle}.

\paragraph{Directions.}
For a small step $h>0$, comparing $F(p+hv)$ with $F(p)$ reveals the sign of the directional derivative $\ip{\nabla F(p)}v$, as long as this derivative is not too close to zero.
Indeed, by smoothness, the finite difference $(F(p+hv)-F(p))/h$ differs from $\ip{\nabla F(p)}v$ by at most $Lh\norm v^2/2$.
The routine $\Dir$ (\Cref{alg:dir} in \Cref{app:direction}) combines such sign tests in three stages following Comparison-GDE of~\citet[Algorithm~2]{zhangli2024}.
First, it determines the sign of every coordinate of the gradient.
Second, it selects a pivot coordinate whose magnitude is close to the largest one, by a tournament in which each comparison tests which of two coordinates is larger in magnitude.
Third, for every other coordinate, it determines the ratio of that coordinate to the pivot coordinate by bisection, where each bisection step is a sign test along a combination of two coordinate directions.
Normalizing the vector of estimated ratios gives the estimate of the gradient direction.
The step sizes of the sign tests are chosen small enough that the finite-difference error does not change any sign that matters for the result.

Formally, let $F:\R^n\to\R$ be differentiable with
\begin{equation}\label{eq:direction-smoothness}
 \norm{\nabla F(u)-\nabla F(v)}\le L\norm{u-v}
 \fa u,v\in\R^n.
\end{equation}
The parameter $\gamma$ of $\Dir(p;L,\gamma,\xi)$ (\Cref{alg:dir}) is an assumed lower bound on $\norm{\nabla F(p)}$ that determines the step sizes of the sign tests.

We use the routine in two contexts: 
with $F=f$ and $n=d$, to estimate the gradient direction at the iterates of \Cref{alg:finite}; 
and with $F$ the restriction of $f$ to a line or a plane, so that $n\in\{1,2\}$, to compute the cuts of the subproblem solver introduced below.

\begin{lemma}[Direction estimation]\label{lem:direction}
Let $F:\R^n\to\R$ be differentiable and satisfy~\eqref{eq:direction-smoothness} for some $L>0$, let $p\in\R^n$, $\gamma>0$, and $\xi\in(0,1]$, and let $\widehat n$ be the output of $\Dir(p;L,\gamma,\xi)$.
Then $\widehat n$ is a unit vector, and the routine uses
\begin{equation}\label{eq:dircost}
 D_n(\xi)=2n-1+(n-1)\left\lceil\log_2\frac{16\sqrt n}{\xi}\right\rceil
 =O\!\left(n\log\frac{2n}{\xi}\right)
\end{equation}
comparisons.
If, in addition, $\norm{\nabla F(p)}\ge\gamma$, then $\norm{\widehat n-\nabla F(p)/\norm{\nabla F(p)}}\le\xi/2$.
\end{lemma}
The proof is given in \Cref{app:direction}.

\paragraph{Small subproblems.}
We solve the ray and plane problems with a cutting-plane method over a ball: 
bisection in dimension one and the central-cut ellipsoid method~\citep{gls} in dimension two.
The routine $\SmallMin$ (\Cref{alg:small} in \Cref{app:small}) replaces the exact gradient that defines each cut by the output of $\Dir$, and it keeps track of the best feasible center by comparisons.

\begin{lemma}[Small subproblems]\label{lem:small}
Let $n\in\{1,2\}$, let $\varphi:\R^n\to\R$ be convex and $L$-smooth, and let $B>0$, $M\ge0$, and $\nu>0$ be such that
\begin{equation}\label{eq:smallbound}
 \norm{\nabla\varphi(u)}\le M\fa u\in C\coloneqq\Ball(0,B).
\end{equation}
Then $\SmallMin(\varphi;B,L,M,\nu)$ returns a point $b\in C$ with $\varphi(b)-\min_{u\in C}\varphi(u)\le\nu$, using $O(1+\log_+(BM/\nu))$ comparisons for $n=1$ and $O([1+\log_+(BM/\nu)]^2)$ comparisons for $n=2$, where $\log_+t\coloneqq\log\max\{1,t\}$.
\end{lemma}

In dimension zero, the feasible set is a single point, which is returned without comparisons.
The proof is given in \Cref{app:small}.

\subsection{The comparison-based method}\label{sec:algorithm}

\Cref{alg:finite} combines the two routines into a method that uses only comparisons.
Each iteration mirrors an iteration of \Cref{alg:ideal}, with the following modifications:
\begin{enumerate}[label=(\roman*),leftmargin=*]
\item \emph{Directions.} The unit gradient $n^k$ is replaced by the estimate $\widehat n^k=\Dir(x^k;L,\gamma,\xi)$, and the memory vector accumulates these estimates.
\item \emph{Inexact searches over balls.} Instead of the whole ray $\{x^k-\alpha\widehat n^k:\alpha\ge0\}$, the ray search considers only the segment $I_k$ of the ray that lies in the ball $\Ball(x^0,B_{k+1})$, and finds a point $z\in I_k$ whose value is within $\rho_{\rm line}$ of the minimum of $f$ over $I_k$.
Likewise, instead of the whole plane $x^0+V_{k+1}$, the plane search considers only the disk $C_{k+1}$ in which this plane meets $\Ball(x^0,B_{k+1})$, and finds a point $w\in C_{k+1}$ whose value is within $\rho_{\rm plane}$ of the minimum of $f$ over $C_{k+1}$.
\item \emph{Keeping the better point.} The method compares $z$ with $x^k$ and keeps the better of the two as $p^{k+1}$; it then compares $w$ with $p^{k+1}$ and keeps the better of the two as $x^{k+1}$. Hence the objective values $f(x^k)$ never increase.
\item \emph{Fixed horizon.} The method runs for $K$ iterations and returns $x^K$.
\end{enumerate}

\begin{algorithm}[tbp]
\caption{Comparison-based normalized subspace method}\label{alg:finite}
\begin{algorithmic}[1]
\Require Comparisons of $f$; $x^0$, $L$, $R$, and $\varepsilon>0$.
\If{$\varepsilon>LR^2$}\label{line:finite-trivial} \State \Return $x^0$. \Comment{$x^0$ is already $\varepsilon$-optimal}\EndIf
\State Set $K=\lceil10\sqrt{LR^2/\varepsilon}\,\rceil$. \Comment{Horizon}\label{line:finite-init}
\State Set $B_k=(16K+3k)R$ for $0\le k\le K$. \Comment{Radii}
\State Set $\gamma=\varepsilon/((19K+1)R)$ and $\xi=1/(64K)$. \Comment{Direction parameters}
\State Set $\rho_{\rm line}=\gamma^2/(8L)$ and $\rho_{\rm plane}=\gamma^2/(5776LK^4)$. \Comment{Search tolerances}\label{line:finite-tol}
\State $m^0=0$.
\For{$k=0,\ldots,K-1$}
\State $\widehat n^k=\Dir(x^k;L,\gamma,\xi)$. \Comment{Estimated unit gradient}\label{line:finite-direction}
\State $m^{k+1}=m^k+\widehat n^k$. \Comment{Memory}
\State $I_k=\{x^k-\alpha\widehat n^k:\alpha\ge0\}\cap\Ball(x^0,B_{k+1})$. \Comment{Segment of the ray}
\State Find $z\in I_k$ with $f(z)\le\min_{I_k}f+\rho_{\rm line}$, using $\SmallMin$.\label{line:finite-ray}
\State Compare $z$ and $x^k$; retain the better point as $p^{k+1}$.\label{line:finite-ray-retain}
\State $V_{k+1}=\spanof\{p^{k+1}-x^0,m^{k+1}\}$. \Comment{Plane directions}\label{line:finite-plane-setup}
\State $C_{k+1}=(x^0+V_{k+1})\cap\Ball(x^0,B_{k+1})$. \Comment{Disk}
\State Find $w\in C_{k+1}$ with $f(w)\le\min_{C_{k+1}}f+\rho_{\rm plane}$, using $\SmallMin$.\label{line:finite-plane-solve}
\State Compare $w$ and $p^{k+1}$; retain the better point as $x^{k+1}$.\label{line:finite-plane-retain}
\EndFor
\State \Return $x^K$.
\end{algorithmic}
\end{algorithm}

\paragraph{Parameters.}
The number of iterations $K=\lceil10\sqrt A\,\rceil$ is chosen so that the bound~\eqref{eq:robust-rate} of \Cref{thm:inexact-certificate} gives $\Delta_K\le\varepsilon$.
The threshold $\gamma$ is a lower bound on the gradient norm at every iterate that is not yet $\varepsilon$-optimal, so that \Cref{lem:direction} applies there; $\xi$ sets the accuracy of the direction estimates.
The radii $B_k$ are of order $KR$, much larger than the distance $R$ from $x^0$ to $x^\star$; their role is explained in the proof sketch below.
The tolerances $\rho_{\rm line}$ and $\rho_{\rm plane}$ are polynomial in $1/K$.

\paragraph{The searches as small subproblems.}
The ray search in line~\ref{line:finite-ray}, which replaces~\eqref{eq:idealray}, and the plane search in line~\ref{line:finite-plane-solve}, which replaces~\eqref{eq:idealplane}, are both carried out by $\SmallMin$: in orthonormal coordinates centered at the midpoint of the interval $I_k$ and at $x^0$, respectively, each search minimizes a convex $L$-smooth function over a ball $\Ball(0,B)$ in dimension $n\le2$ with $B\le B_{k+1}$, and the gradient bound $M=L(B_{k+1}+R)$ of \Cref{lem:small} holds, since every point of $\Ball(x^0,B_{k+1})$ lies within distance $B_{k+1}+R$ of $x^\star$, where the gradient vanishes (see \Cref{app:outer} for the coordinates).

\begin{theorem}[Comparison-based method]\label{thm:finite}
Under \Cref{ass:setup}, \Cref{alg:finite} returns a point satisfying~\eqref{eq:target}.
For $\varepsilon>LR^2$, it makes no comparisons;
for $0<\varepsilon\le LR^2$, it uses
\begin{equation}\label{eq:cost}
 O\!\left(K[d\log(dK)+\log^2K]\right)
 =O\!\left(\sqrt A\,[d\log(2dA)+\log^2(2A)]\right)
 =\widetilde O\!\left(d\sqrt{A}\right)
\end{equation}
comparisons, where $K=\lceil10\sqrt A\,\rceil$.
\end{theorem}

\begin{proofsketch}
The full proof is given in \Cref{app:outer}.
It verifies the conditions of \Cref{thm:inexact-certificate} for the iterates of \Cref{alg:finite} with $u^k=\widehat n^k$ and
\begin{equation}\label{eq:finite-certificate-budgets}
 \vartheta_k=\frac{9k}{64K},\qquad r_k=\frac RK,\qquad
 \delta_k=\frac1{64K},\qquad \eta=\frac1{4L}\fa 0\le k<K.
\end{equation}
Since the objective values never increase, the conditions are only needed at iterations with $\Delta_k>\varepsilon$.

\emph{Direction condition.}
If $\Delta_k>\varepsilon$, convexity gives $\varepsilon<\Delta_k\le G_k\norm{x^k-x^\star}\le G_k(B_k+R)$, hence $G_k>\gamma$.
\Cref{lem:direction} then bounds the direction error by $\xi/2\le\delta_k$.

\emph{Memory and anchor conditions.}
For $k\ge1$, the point $x^k$ is within $\rho_{\rm plane}$ of the minimum of $f$ over the disk $C_k$, but its gradient need not be orthogonal to the plane, as it would be for an exact minimizer.
Instead, near-optimality bounds the Frank--Wolfe gap $\max_{u\in C_k}\ip{g_k}{x^k-u}$ by $\gamma R/K$ (\Cref{lem:vi}).
Taking $u=x^0$ gives the anchor condition.
If $m^k\ne0$, taking $u=x^0\pm B_k m^k/\norm{m^k}$, which lie in $C_k$, gives $|\ip{g_k}{m^k}|\le(\gamma R/K+G_kR)\norm{m^k}/B_k$.
Dividing by $G_k>\gamma$, accounting for the direction error, and using $B_k\ge16KR$ and $\norm{m^k}\le k$ shows that $\ip{\widehat n^k}{m^k}\le9k/(64K)$, which is the memory condition.
For $k=0$, where $x^0-x^0=0$ and $m^0=0$, both conditions hold trivially, and for $k\ge1$ with $m^k=0$, so does the memory condition.
This is where the large radii are needed: the larger the disk, the closer near-optimality comes to orthogonality.

\emph{Descent condition.}
Since the radii grow by $3R$ per iteration, the point $x^k-(G_k/L)\widehat n^k$ lies on the segment $I_k$.
By smoothness, this point decreases $f$ by nearly $G_k^2/(2L)$, and the ray search loses at most $\rho_{\rm line}\le G_k^2/(8L)$ of this decrease; the plane search can only decrease $f$ further.
This leaves a decrease of at least $G_k^2/(4L)$.

\emph{Conclusion.}
With the values~\eqref{eq:finite-certificate-budgets}, we have $S_K\le3R\sqrt K$, and \Cref{thm:inexact-certificate} gives $\Delta_K\le\max\{\varepsilon,\allowbreak 36LR^2/K^2\}$, which is $\varepsilon$ since $K\ge10\sqrt A$.
Each iteration uses $O(d\log(dK))$ comparisons for the direction (\Cref{lem:direction}), and $O(\log^2K)$ comparisons for the two searches and the two retention steps (\Cref{lem:small}), because $BM$ divided by either tolerance is polynomial in $K$.
\end{proofsketch}

\paragraph{Near-optimality.}
Since two function values determine one comparison, every lower bound on the number of exact function values needed to reach accuracy $\varepsilon$ is also a lower bound on the number of comparisons, up to a factor of two.
Applied to the lower bound of~\citet[Theorem~2.1]{dorn2026}, this gives the following corollary.

\begin{corollary}[Lower bound for comparison methods]\label{cor:lower-bound}
Let $d\ge1$, $L,R>0$, and $0<\varepsilon\le LR^2/32$, that is, $A\ge32$.
Suppose that a method, possibly randomized, accesses $f$ only through the comparison oracle~\eqref{eq:oracle}, receives $x^0$, $L$, $R$, and $\varepsilon$, makes at most $T$ comparisons in every run, and, for every $x^0\in\R^d$ and every $f$ satisfying \Cref{ass:setup}, returns a point satisfying~\eqref{eq:target} with probability at least $2/3$.
The method uses a random seed, independent of $f$, on a standard Borel probability space (a constant seed for deterministic methods), and its queries and its output are Borel functions of the seed and the previous replies, defined for every finite sequence of replies.
Then
\begin{equation}\label{eq:lower-bound}
 T\ge\frac{d}{8192}\min\bigl\{d,\sqrt A\bigr\}.
\end{equation}
\end{corollary}

\begin{proof}
\citet[Section~2]{dorn2026} consider methods that query exact function values, possibly using an independent random seed, and the class of convex, $L$-smooth functions with a minimizer in $\Ball(0,R)$, that is, the functions satisfying \Cref{ass:setup} with $x^0=0$.
From the comparison method, run with $x^0=0$ and the same seed, we obtain such a method:
it answers each comparison $O_f(u,v)$ by querying $f(u)$ and $f(v)$, and it returns $0$ if the comparison method has not returned after $T$ comparisons.
This method uses at most $2T$ function values for every sequence of function values, and it has the same success probability as the comparison method for every function in the class.
Since the queries and the output of the comparison method are Borel functions of its seed and its replies, defined for every sequence of replies, and each reply is a Borel function of two function values, the queries and the output of the simulation are Borel functions of the seed and the function values, defined for every sequence of function values, as required by~\citet[Section~2]{dorn2026}.
Therefore, \citet[Theorem~2.1]{dorn2026}, with their constant $C_\star=4096$, give $2T\ge\frac d{4096}\min\{d,\sqrt A\}$, which is~\eqref{eq:lower-bound}.
\end{proof}

Hence, for $32\le A\le d^2$, the bound of \Cref{thm:finite} is optimal up to logarithmic factors, since every such comparison method needs $\Omega(d\sqrt A)$ comparisons by~\eqref{eq:lower-bound}.
More precisely, $d\ge6$ in this regime, and $\log^2(2A)\le2d\log(2dA)$ and $\log(2dA)\le3\log(2d)$, so the bound of \Cref{thm:finite} is $O(d\sqrt A\log d)$, within a factor $O(\log d)$ of~\eqref{eq:lower-bound}.
For $A>d^2$, the lower bound becomes $\Omega(d^2)$, and the cutting-plane method of~\citet[Theorem~3]{zhangli2024}, which uses $\widetilde O(d^2)$ comparisons, is optimal up to logarithmic factors. 
\Cref{prop:cutting} in \Cref{app:cutting} gives a self-contained version with $O(d^2\log A\log(dA))$ comparisons.
Since $d$ and $A$ are known, one can run \Cref{alg:finite} if $A\le d^2$ and the cutting-plane method otherwise. 
Together with \Cref{cor:lower-bound}, the comparison complexity of the class is therefore $\widetilde\Theta(d\min\{d,\sqrt A\})$ for $A\ge32$.

\paragraph{Discussion.}\label{par:discussion}
The guarantee of \Cref{thm:finite} uses a prescribed horizon: \Cref{alg:finite} runs for $K=\lceil10\sqrt A\,\rceil$ iterations and has no stopping certificate.
Once the gap drops below $\varepsilon$, monotonicity keeps it there, although the method does not observe this.
This number of iterations is dictated by the worst case of the analysis, and it is about seven times the $\sqrt{2A}$ iterations that suffice for the exact normalized subspace method (\Cref{thm:ideal}).
Each iteration may be expensive as the plane search must reach a tolerance of order $LR^2/K^{10}$.

The analysis assumes exact comparisons and exact real arithmetic.
The inner tolerances are tiny (for $A=10^2$, the plane tolerance is about $5\cdot10^{-23}LR^2$), so straightforward floating-point evaluation of $f$ may fail to resolve the required comparisons, depending on the scale and offset of $f$; 
a numerical implementation may need higher precision or a specialized comparison routine, and \Cref{thm:finite} gives no finite-precision guarantee.

\paragraph{Noisy comparisons.}
\Cref{alg:finite} is well defined and makes at most $N=O(K[d\log(dK)+\log^2K])$ comparisons for every sequence of replies, correct or not (\Cref{lem:noisy-budget}).
If each reply is correct with probability at least a known $p>1/2$ given all previous queries and replies~\eqref{eq:noisy-model}, it therefore suffices to repeat every comparison $O((p-1/2)^{-2}\log(N/\delta))$ times and to take majority votes: 
with probability at least $1-\delta$, all majority votes are correct, and \Cref{thm:finite} applies (\Cref{cor:noisy}).
In this model, this improves the dependence of the bound of~\citet{saha} on the accuracy from $1/\varepsilon$ to $1/\sqrt\varepsilon$.
If the probability of a correct reply approaches $1/2$ as the two function values approach each other, then no $p>1/2$ as in~\eqref{eq:noisy-model} exists and the repetition argument does not apply.

\section{Conclusion}\label{sec:conclusion}

A comparison oracle reveals neither function values nor gradient norms, yet it suffices for accelerated smooth convex optimization.
We gave a deterministic comparison-based method that finds an $\varepsilon$-optimal point with $\widetilde O(d\sqrt A)$ comparisons, which matches the accelerated $\sqrt A$ dependence of methods with exact function values and is optimal up to logarithmic factors when $32\le A\le d^2$.
The key is that the analysis of the exact normalized subspace method uses only gradient directions, one- and two-dimensional minimizations, and relaxed versions of its orthogonality and descent properties, all of which comparisons can provide.
The method needs the constants $L$ and $R$ and the target accuracy in advance, its guarantees assume exact arithmetic, and its search tolerances decrease polynomially in $1/K$, which may make a floating-point implementation demanding; noisy comparisons are covered when each reply is correct with probability at least a known $p>1/2$ given the previous replies (\Cref{cor:noisy}).
Removing these requirements, giving finite-precision guarantees, closing the logarithmic gap to the lower bound, and determining the optimal comparison complexity under persistent errors~\citep{braverman2008,lobanov}, or under noise whose probability of a correct reply approaches $1/2$ as the function values approach each other, remain open.

\section*{Acknowledgments}
The research reported in this publication was supported by funding from King Abdullah University of Science and Technology (KAUST): i) KAUST Baseline Research Scheme, ii) CRG Grant ORFS-CRG12-2024-6460, and iii) Center of Excellence for Generative AI, under award number 5940.
This paper was developed in collaboration with AI.

{
\small
\bibliographystyle{plainnat}
\bibliography{references}

\begin{thebibliography}{33}
\providecommand{\natexlab}[1]{#1}
\providecommand{\url}[1]{\texttt{#1}}
\expandafter\ifx\csname urlstyle\endcsname\relax
  \providecommand{\doi}[1]{doi: #1}\else
  \providecommand{\doi}{doi: \begingroup \urlstyle{rm}\Url}\fi

\bibitem[Allen-Zhu and Hazan(2016)]{allenzhu2016}
Zeyuan Allen-Zhu and Elad Hazan.
\newblock Optimal black-box reductions between optimization objectives.
\newblock In \emph{Advances in Neural Information Processing Systems},
  volume~29, pages 1614--1622, 2016.

\bibitem[Bergou et~al.(2020)Bergou, Gorbunov, and Richt{\'a}rik]{bergou2020}
El~Houcine Bergou, Eduard Gorbunov, and Peter Richt{\'a}rik.
\newblock Stochastic three points method for unconstrained smooth minimization.
\newblock \emph{SIAM Journal on Optimization}, 30\penalty0 (4):\penalty0
  2726--2749, 2020.

\bibitem[Bradley and Terry(1952)]{bradleyterry1952}
Ralph~Allan Bradley and Milton~E. Terry.
\newblock Rank analysis of incomplete block designs: {I}. {T}he method of
  paired comparisons.
\newblock \emph{Biometrika}, 39\penalty0 (3/4):\penalty0 324--345, 1952.

\bibitem[Braverman and Mossel(2008)]{braverman2008}
Mark Braverman and Elchanan Mossel.
\newblock Noisy sorting without resampling.
\newblock In \emph{Proceedings of the Nineteenth Annual ACM-SIAM Symposium on
  Discrete Algorithms}, pages 268--276, 2008.

\bibitem[Bubeck(2015)]{bubeck2015}
S{\'e}bastien Bubeck.
\newblock Convex optimization: Algorithms and complexity.
\newblock \emph{Foundations and Trends in Machine Learning}, 8\penalty0
  (3--4):\penalty0 231--358, 2015.

\bibitem[Cai et~al.(2022)Cai, McKenzie, Yin, and Zhang]{cai}
HanQin Cai, Daniel McKenzie, Wotao Yin, and Zhenliang Zhang.
\newblock A one-bit, comparison-based gradient estimator.
\newblock \emph{Applied and Computational Harmonic Analysis}, 60:\penalty0
  242--266, 2022.

\bibitem[Dorn et~al.(2026)Dorn, Dvinskikh, Loginov, Osmanov, Shestakov,
  Tupitsa, and Gasnikov]{dorn2026}
Yuriy Dorn, Darina Dvinskikh, Timofei Loginov, Osman Osmanov, Aleksandr
  Shestakov, Nazarii Tupitsa, and Alexander Gasnikov.
\newblock Lower bounds for gradient-free convex optimization and convex-concave
  saddle-point problems, 2026.
\newblock arXiv:2610.00345v1.

\bibitem[Feige et~al.(1994)Feige, Raghavan, Peleg, and Upfal]{feige1994}
Uriel Feige, Prabhakar Raghavan, David Peleg, and Eli Upfal.
\newblock Computing with noisy information.
\newblock \emph{SIAM Journal on Computing}, 23\penalty0 (5):\penalty0
  1001--1018, 1994.

\bibitem[Golovin et~al.(2020)Golovin, Karro, Kochanski, Lee, Song, and
  Zhang]{golovin2020}
Daniel Golovin, John Karro, Greg Kochanski, Chansoo Lee, Xingyou Song, and
  Qiuyi Zhang.
\newblock Gradientless descent: High-dimensional zeroth-order optimization.
\newblock In \emph{International Conference on Learning Representations}, 2020.

\bibitem[Gr{\"o}tschel et~al.(1988)Gr{\"o}tschel, Lov{\'a}sz, and
  Schrijver]{gls}
Martin Gr{\"o}tschel, L{\'a}szl{\'o} Lov{\'a}sz, and Alexander Schrijver.
\newblock \emph{Geometric Algorithms and Combinatorial Optimization}.
\newblock Springer, Berlin, 1988.

\bibitem[Hoeffding(1963)]{hoeffding1963}
Wassily Hoeffding.
\newblock Probability inequalities for sums of bounded random variables.
\newblock \emph{Journal of the American Statistical Association}, 58\penalty0
  (301):\penalty0 13--30, 1963.

\bibitem[Jaggi(2013)]{jaggi2013}
Martin Jaggi.
\newblock Revisiting {Frank--Wolfe}: Projection-free sparse convex
  optimization.
\newblock In \emph{Proceedings of the 30th International Conference on Machine
  Learning}, volume~28 of \emph{Proceedings of Machine Learning Research},
  pages 427--435, 2013.

\bibitem[Jamieson et~al.(2012)Jamieson, Nowak, and Recht]{jamieson2012}
Kevin~G. Jamieson, Robert~D. Nowak, and Benjamin Recht.
\newblock Query complexity of derivative-free optimization.
\newblock In \emph{Advances in Neural Information Processing Systems},
  volume~25, pages 2672--2680, 2012.

\bibitem[Karabag et~al.(2021)Karabag, Neary, and Topcu]{karabag}
Mustafa~O. Karabag, Cyrus Neary, and Ufuk Topcu.
\newblock Smooth convex optimization using sub-zeroth-order oracles.
\newblock In \emph{Proceedings of the AAAI Conference on Artificial
  Intelligence}, volume~35, pages 3815--3822, 2021.

\bibitem[Kumagai(2017)]{kumagai2017}
Wataru Kumagai.
\newblock Regret analysis for continuous dueling bandit.
\newblock In \emph{Advances in Neural Information Processing Systems},
  volume~30, pages 1489--1498, 2017.

\bibitem[Kuruzov and Stonyakin(2026)]{kuruzov}
Ilya Kuruzov and Fedor Stonyakin.
\newblock On some versions of subspace optimization methods with inexact
  gradient information, 2026.
\newblock arXiv:2407.01051v3.

\bibitem[Levy(2017)]{levy2017}
Kfir~Y. Levy.
\newblock Online to offline conversions, universality and adaptive minibatch
  sizes.
\newblock In \emph{Advances in Neural Information Processing Systems},
  volume~30, pages 1613--1622, 2017.

\bibitem[Lobanov et~al.(2024)Lobanov, Gasnikov, and Krasnov]{lobanov}
Aleksandr Lobanov, Alexander Gasnikov, and Andrei Krasnov.
\newblock Acceleration exists! {O}ptimization problems when oracle can only
  compare objective function values.
\newblock In \emph{Advances in Neural Information Processing Systems},
  volume~37, pages 21428--21460, 2024.

\bibitem[Narkiss and Zibulevsky(2005)]{narkiss}
Guy Narkiss and Michael Zibulevsky.
\newblock Sequential subspace optimization method for large-scale unconstrained
  problems.
\newblock Technical report, Technion -- Israel Institute of Technology, 2005.

\bibitem[Nemirovski(1982)]{nemirovski}
Arkadi Nemirovski.
\newblock {ORTH}-method for smooth convex optimization.
\newblock \emph{Izvestiya Akademii Nauk SSSR, Tekhnicheskaya Kibernetika},
  \penalty0 (2), 1982.
\newblock In Russian.

\bibitem[Nemirovsky and Yudin(1983)]{nemirovskyyudin1983}
Arkadi~S. Nemirovsky and David~B. Yudin.
\newblock \emph{Problem Complexity and Method Efficiency in Optimization}.
\newblock Wiley, Chichester, 1983.
\newblock Translated from the Russian edition (Nauka, 1979).

\bibitem[Nesterov(1983)]{nesterov1983}
Yurii Nesterov.
\newblock A method of solving a convex programming problem with convergence
  rate {$O(1/k^2)$}.
\newblock \emph{Doklady Akademii Nauk SSSR}, 269\penalty0 (3):\penalty0
  543--547, 1983.
\newblock In Russian.

\bibitem[Nesterov and Spokoiny(2017)]{nesterovspokoiny2017}
Yurii Nesterov and Vladimir Spokoiny.
\newblock Random gradient-free minimization of convex functions.
\newblock \emph{Foundations of Computational Mathematics}, 17\penalty0
  (2):\penalty0 527--566, 2017.

\bibitem[Rockafellar(1970)]{rockafellar1970}
R.~Tyrrell Rockafellar.
\newblock \emph{Convex Analysis}.
\newblock Princeton University Press, Princeton, NJ, 1970.

\bibitem[Saha et~al.(2021)Saha, Koren, and Mansour]{saha}
Aadirupa Saha, Tomer Koren, and Yishay Mansour.
\newblock Dueling convex optimization.
\newblock In \emph{Proceedings of the 38th International Conference on Machine
  Learning}, volume 139 of \emph{Proceedings of Machine Learning Research},
  pages 9245--9254, 2021.

\bibitem[Scheinberg and Xiong(2026)]{scheinberg}
Katya Scheinberg and Zikai Xiong.
\newblock Function-free optimization via comparison oracles, 2026.
\newblock arXiv:2604.26867v2.

\bibitem[Stich et~al.(2013)Stich, M{\"u}ller, and G{\"a}rtner]{stich2013}
Sebastian~U. Stich, Christian~L. M{\"u}ller, and Bernd G{\"a}rtner.
\newblock Optimization of convex functions with random pursuit.
\newblock \emph{SIAM Journal on Optimization}, 23\penalty0 (2):\penalty0
  1284--1309, 2013.

\bibitem[Tao et~al.(2026)Tao, Zhang, Wang, Zhang, and Li]{tao}
Xiwen Tao, Chenyi Zhang, Helin Wang, Yexin Zhang, and Tongyang Li.
\newblock Gradient testing and estimation by comparisons, 2026.
\newblock arXiv:2405.11454v3.

\bibitem[Thurstone(1927)]{thurstone1927}
Louis~L. Thurstone.
\newblock A law of comparative judgment.
\newblock \emph{Psychological Review}, 34\penalty0 (4):\penalty0 273--286,
  1927.

\bibitem[Tyurin et~al.(2026{\natexlab{a}})Tyurin, Pham, Khrulkov, and
  Oseledets]{longtails}
Alexander Tyurin, QuocViet Pham, Valentin Khrulkov, and Ivan Oseledets.
\newblock \mbox{LongTails}: Automated research in optimization.
\newblock GitHub repository, commit 7c8f00f, 2026{\natexlab{a}}.
\newblock URL \url{https://github.com/k3nfalt/longtails}.

\bibitem[Tyurin et~al.(2026{\natexlab{b}})Tyurin, Pham, Khrulkov, and
  Oseledets]{tyurin2026longtails}
Alexander Tyurin, QuocViet Pham, Valentin Khrulkov, and Ivan Oseledets.
\newblock {LongTails}: Automated research in optimization, 2026{\natexlab{b}}.
\newblock arXiv:2611.13872v1.

\bibitem[Wu et~al.(2026)Wu, Zhang, Zhang, Li, Zheng, Fang, Li, and Lin]{wu2026}
Wendao Wu, Haihan Zhang, Chenheng Zhang, Yanyi Li, Chunyuan Zheng, Cong Fang,
  Haoxuan Li, and Zhouchen Lin.
\newblock Near-optimal deterministic exact-value complexity for smooth convex
  optimization, 2026.
\newblock arXiv:2609.18230v2.

\bibitem[Zhang and Li(2024)]{zhangli2024}
Chenyi Zhang and Tongyang Li.
\newblock Comparisons are all you need for optimizing smooth functions, 2024.
\newblock arXiv:2405.11454v1.

\end{thebibliography}
}

\appendix
\section{Proof of \texorpdfstring{\Cref{thm:inexact-certificate}}{the theorem on relaxed conditions}}\label{app:ideal}

The proof of \Cref{thm:inexact-certificate} follows the proof of \Cref{thm:ideal}, with each exact identity replaced by the corresponding inequality in~\eqref{eq:robust-conditions}.

\begin{proof}[Proof of \Cref{thm:inexact-certificate}]
If $\Delta_j\le\varepsilon$ for some $0\le j\le K$, then $\Delta_K\le\Delta_j\le\varepsilon$ by monotonicity, and~\eqref{eq:robust-rate} holds.
We may therefore assume that $\Delta_k>\varepsilon\ge0$ for all $0\le k\le K$.
Then the conditions~\eqref{eq:robust-conditions} hold for every $0\le k<K$, and all gaps $\Delta_k$ and gradient norms $G_k$ are positive.

\pstep{Norm of the memory vector}
Since $u^k$ is a unit vector, the definition of $m^{k+1}$ and~\eqref{eq:cond-memory} give $\norm{m^{k+1}}^2-\norm{m^k}^2=2\ip{u^k}{m^k}+\norm{u^k}^2\le1+2\vartheta_k$.
Summing over $0\le k<K$ and using $m^0=0$ gives
\begin{equation}\label{eq:robust-memory-bound}
 \norm{m^K}^2\le K+2\sum_{k=0}^{K-1}\vartheta_k .
\end{equation}

\pstep{Upper bound on normalized gaps}
By convexity, $\Delta_k\le\ip{g_k}{x^k-x^\star}$.
We split $x^k-x^\star=(x^k-x^0)+(x^0-x^\star)$, divide by $G_k$, and use~\eqref{eq:cond-direction}, \eqref{eq:cond-anchor}, and the Cauchy--Schwarz inequality:
\begin{equation}\label{eq:robust-one-gap}
 \frac{\Delta_k}{G_k}
 \le\frac{\ip{g_k}{x^k-x^0}}{G_k}+\ip{u^k}{x^0-x^\star}+\ip{\frac{g_k}{G_k}-u^k}{x^0-x^\star}
 \le r_k+\ip{u^k}{x^0-x^\star}+R\delta_k .
\end{equation}
Summing over $k$, and using $m^K=\sum_{k<K}u^k$, the Cauchy--Schwarz inequality, and~\eqref{eq:robust-memory-bound},
\begin{equation}\label{eq:robust-gap-sum}
 \sum_{k=0}^{K-1}\frac{\Delta_k}{G_k}
 \le\ip{m^K}{x^0-x^\star}+\sum_{k=0}^{K-1}(r_k+R\delta_k)
 \le R\norm{m^K}+\sum_{k=0}^{K-1}(r_k+R\delta_k)
 \le S_K .
\end{equation}

\pstep{Lower bound on normalized gaps}
Define $D_k\coloneqq\Delta_k-\Delta_{k+1}$ and $t_k\coloneqq D_k/\Delta_k^2$.
By~\eqref{eq:cond-direction} and~\eqref{eq:cond-descent}, $D_k\ge\eta G_k^2>0$, so $t_k>0$ and $\Delta_k/G_k\ge\sqrt\eta\,\Delta_k/\sqrt{D_k}=\sqrt\eta\,t_k^{-1/2}$ for $0\le k<K$.
Since $0<\Delta_{k+1}\le\Delta_k$, the telescoping bound $\sum_{k<K}t_k\le1/\Delta_K$ and Jensen's inequality follow as in the proof of \Cref{thm:ideal}, and they give $\sum_{k<K}\Delta_k/G_k\ge\sqrt\eta\,K^{3/2}\sqrt{\Delta_K}$.

\pstep{Conclusion}
Combining the upper and lower bounds gives $\sqrt\eta\,K^{3/2}\sqrt{\Delta_K}\le S_K$, that is, $\Delta_K\le S_K^2/(\eta K^3)$, which proves~\eqref{eq:robust-rate}.
\end{proof}

\section{Proof of \texorpdfstring{\Cref{thm:finite}}{the comparison-based method}}\label{app:outer}

We first recall two standard facts about convex functions with Lipschitz gradient, which are used in the remaining proofs; 
we include their proofs for completeness.

\begin{lemma}[Smoothness facts]\label{lem:smooth}
Under \Cref{ass:setup}, for every $x,s\in\R^d$,
\begin{equation}\label{eq:upper}
 \left|f(x+s)-f(x)-\ip{\nabla f(x)}s\right|\le\frac L2\norm{s}^2.
\end{equation}
Moreover, for every $x\in\R^d$,
\begin{equation}\label{eq:nonexp}
 \norm{x-\nabla f(x)/L-x^\star}^2
 \le\norm{x-x^\star}^2-\frac{\norm{\nabla f(x)}^2}{L^2}.
\end{equation}
\end{lemma}

Inequality~\eqref{eq:upper} only uses that the gradient is $L$-Lipschitz.
We therefore also apply it to other functions with this property, such as the function $F$ in~\eqref{eq:direction-smoothness} and the restrictions of $f$ to lines and planes.

\begin{proof}
\pstep{Upper model}
By the fundamental theorem of calculus and the Lipschitz continuity of the gradient, $|f(x+s)-f(x)-\ip{\nabla f(x)}s|\le\int_0^1\norm{\nabla f(x+ts)-\nabla f(x)}\norm{s}\,dt\le L\norm{s}^2\int_0^1t\,dt=\frac L2\norm{s}^2$.

\pstep{Co-coercivity}
Let $x,y\in\R^d$.
The function $h(z)\coloneqq f(z)-\ip{\nabla f(y)}z$ is convex and $L$-smooth, and $\nabla h(y)=0$, so $y$ is a global minimizer of $h$.
Applying~\eqref{eq:upper} to $h$ at $x$ with $s=-\nabla h(x)/L$ gives $h(y)\le h(x-\nabla h(x)/L)\le h(x)-\norm{\nabla h(x)}^2/L+\norm{\nabla h(x)}^2/(2L)=h(x)-\norm{\nabla h(x)}^2/(2L)$.
Since $h(x)-h(y)=f(x)-f(y)-\ip{\nabla f(y)}{x-y}$ and $\nabla h(x)=\nabla f(x)-\nabla f(y)$, this reads $f(x)-f(y)-\ip{\nabla f(y)}{x-y}\ge\norm{\nabla f(x)-\nabla f(y)}^2/(2L)$.
Adding this inequality to the same inequality with $x$ and $y$ exchanged, the function values cancel, and we obtain
\begin{equation}\label{eq:coco}
 \ip{\nabla f(x)-\nabla f(y)}{x-y}
 \ge\frac{\norm{\nabla f(x)-\nabla f(y)}^2}{L}.
\end{equation}

\pstep{Gradient step}
We apply~\eqref{eq:coco} with $y=x^\star$, where $\nabla f(x^\star)=0$.
Expanding the square gives $\norm{x-\nabla f(x)/L-x^\star}^2=\norm{x-x^\star}^2-\frac2L\ip{\nabla f(x)}{x-x^\star}+\norm{\nabla f(x)}^2/L^2\le\norm{x-x^\star}^2-\norm{\nabla f(x)}^2/L^2$, which is~\eqref{eq:nonexp}.
\end{proof}

\paragraph{Parameters.}
If $\varepsilon>LR^2$, \Cref{alg:finite} returns $x^0$, which is $\varepsilon$-optimal, since by~\eqref{eq:upper} and \Cref{ass:setup}, $f(x^0)-f^\star\le\frac L2\norm{x^0-x^\star}^2\le LR^2/2$.
For $0<\varepsilon\le LR^2$, lines~\ref{line:finite-init}--\ref{line:finite-tol} of \Cref{alg:finite} set the parameters below; for the analysis, we also name $A$, $B_{\max}$, $q$, and $\Gamma$:
\begin{equation}\label{eq:params}
\begin{gathered}
 A=\frac{LR^2}{\varepsilon},\quad K=\lceil10\sqrt A\rceil,\quad
 B_k=(16K+3k)R\fa 0\le k\le K,\\
 B_{\max}=19KR,\quad q=B_{\max}+R,\quad
 \gamma=\frac{\varepsilon}{q},\quad\xi=\frac{1}{64K},\quad\Gamma=\frac{\gamma R}{K},\\
 \rho_{\rm line}=\frac{\gamma^2}{8L},\qquad
 \rho_{\rm plane}=\frac{\Gamma^2}{16LB_{\max}^2}=\frac{\gamma^2}{5776LK^4}.
\end{gathered}
\end{equation}
Since $\varepsilon\le LR^2$ and $q\ge19KR$, we have $\gamma\le LR/(19K)$, hence $\Gamma\le4LB_{\max}^2$, and therefore also $\rho_{\rm plane}\le\Gamma/4$.

\paragraph{Implementation of the searches.}
Iteration $k$ of \Cref{alg:finite} only searches the ray inside the ball $\Ball(x^0,B_{k+1})$.
With $s_k\coloneqq x^k-x^0$, the segment $I_k$ consists of the points $x^k-\alpha\widehat n^k$ with $0\le\alpha\le H_k$, where
\begin{equation}\label{eq:segmentlength}
 c_k=\ip{s_k}{\widehat n^k},\qquad
 H_k=c_k+\sqrt{c_k^2+B_{k+1}^2-\norm{s_k}^2}.
\end{equation}
We center this segment by setting
\begin{equation}\label{eq:rayrestriction}
 \psi_k(u)=f\!\left(x^k-\left(\frac{H_k}{2}+u\right)\widehat n^k\right)
 \fa u\in\R,
\end{equation}
so that the segment corresponds to $|u|\le H_k/2$.
The plane basis $U$ is obtained by applying the Gram--Schmidt procedure first to $p^{k+1}-x^0$ and then to $m^{k+1}$, discarding zero residuals.
Its columns are orthonormal and span $V_{k+1}$, which has dimension zero, one, or two.
With $M_k\coloneqq L(B_{k+1}+R)$, line~\ref{line:finite-ray} of \Cref{alg:finite} computes $u=\SmallMin(\psi_k;H_k/2,L,M_k,\rho_{\rm line})$ and sets $z=x^k-(H_k/2+u)\widehat n^k$, and line~\ref{line:finite-plane-solve} computes $y=\SmallMin(\varphi_k;B_{k+1},L,M_k,\rho_{\rm plane})$ for $\varphi_k(y)\coloneqq f(x^0+Uy)$ and sets $w=x^0+Uy$.
\Cref{lem:ray} and~\eqref{eq:restrictedbound} below show that these calls satisfy the assumptions of \Cref{lem:small}.
All parameters in~\eqref{eq:params} are known and positive.

The plane needs a much smaller tolerance than the line.
On the plane, objective suboptimality is converted into a bound on a Frank--Wolfe gap, which loses a square (\Cref{lem:vi}), while the line search only needs to preserve the decrease.
For $K$ of order $\sqrt A$, the two tolerances are of orders $LR^2/K^6$ and $LR^2/K^{10}$, respectively; they enter the number of comparisons only through logarithms.

\begin{lemma}[Objective suboptimality controls the Frank--Wolfe gap]\label{lem:vi}
Under \Cref{ass:setup}, let $C=(x^0+V)\cap\Ball(x^0,B)$ for a linear subspace $V$ and $B>0$, and let $\nu\ge0$.
If $p\in C$ satisfies $f(p)\le\min_{u\in C}f(u)+\nu$, then its Frank--Wolfe gap~\citep[Eq.~(2)]{jaggi2013} satisfies
\begin{equation}\label{eq:fwgap}
 \omega\coloneqq\max_{u\in C}\ip{\nabla f(p)}{p-u}\ge0,
 \qquad \min\left\{\frac{\omega^2}{8LB^2},\frac{\omega}{2}\right\}\le\nu.
\end{equation}
If, in addition, the parameters are as in~\eqref{eq:params}, $\nu\le\rho_{\rm plane}$, and $B\le B_{\max}$, then
\begin{equation}\label{eq:vi}
 \ip{\nabla f(p)}{p-u}\le\Gamma\fa u\in C.
\end{equation}
\end{lemma}

\begin{proof}
Since $C$ is compact and contains $p$, the maximum in~\eqref{eq:fwgap} is attained at some $u\in C$, and $\omega\ge0$ (take $u=p$).
For $\alpha\in[0,1]$, the point $p+\alpha(u-p)$ lies in $C$, and $\norm{u-p}\le2B$.
By~\eqref{eq:upper}, $f(p+\alpha(u-p))\le f(p)-\alpha\omega+2LB^2\alpha^2$.
If $\omega=0$, there is nothing to prove.
If $0<\omega\le4LB^2$, the choice $\alpha=\omega/(4LB^2)\in(0,1]$ gives a decrease $f(p)-f(p+\alpha(u-p))\ge\omega^2/(8LB^2)$.
If $\omega>4LB^2$, the choice $\alpha=1$ gives a decrease $f(p)-f(u)\ge\omega-2LB^2>\omega/2$.
In both cases, the decrease is at most $f(p)-\min_Cf\le\nu$ by assumption, which proves~\eqref{eq:fwgap}.

For the second claim, suppose that $\omega>\Gamma$.
By~\eqref{eq:params}, $\rho_{\rm plane}\le\Gamma^2/(16LB_{\max}^2)$ and $\rho_{\rm plane}\le\Gamma/4$, so
\begin{equation}\label{eq:vicontra}
 \frac{\omega^2}{8LB^2}>\frac{\Gamma^2}{8LB_{\max}^2}\ge2\rho_{\rm plane}\ge2\nu,
 \qquad \frac{\omega}{2}>\frac{\Gamma}{2}\ge2\rho_{\rm plane}\ge2\nu .
\end{equation}
Thus the minimum in~\eqref{eq:fwgap} exceeds $2\nu\ge\nu$ (and is positive if $\nu=0$), a contradiction.
Hence $\omega\le\Gamma$, which is~\eqref{eq:vi}.
\end{proof}

\begin{lemma}[The ray search is a centered one-dimensional problem]\label{lem:ray}
Under \Cref{ass:setup} and~\eqref{eq:params}, let $\norm{x^k-x^0}\le B_k$, and let $\widehat n^k$ be a unit vector.
Then the intersection of the ray $\{x^k-\alpha\widehat n^k:\alpha\ge0\}$ with $\Ball(x^0,B_{k+1})$ is the image of the interval $|u|\le H_k/2$ under the affine map in~\eqref{eq:rayrestriction}, and $0<H_k\le2B_{k+1}$.
The function $\psi_k$ is convex and $L$-smooth on $\R$, and $|\psi_k'(u)|\le L(B_{k+1}+R)$ for $|u|\le H_k/2$.
\end{lemma}

\begin{proof}
We write $s=s_k$, $c=c_k$, $B=B_{k+1}$, and $v=\widehat n^k$.
Expanding the square, $\norm{s-\alpha v}^2\le B^2$ if and only if $\alpha^2-2c\alpha+\norm{s}^2-B^2\le0$.
Since $\norm{s}\le B_k<B$, the product of the two roots $c\pm\sqrt{c^2+B^2-\norm{s}^2}$ is negative, so the roots have opposite signs.
Intersecting the interval between them with $\alpha\ge0$ gives $[0,H_k]$ with $H_k>0$, by~\eqref{eq:segmentlength}.
Both endpoints of the segment lie in the ball, so its length $H_k$ is at most the diameter $2B$.
The substitution $\alpha=H_k/2+u$ gives~\eqref{eq:rayrestriction}.

Composition with the unit-speed affine map $u\mapsto x^k-(H_k/2+u)v$ preserves convexity and $L$-smoothness.
At a point $z$ of the segment, the chain rule, $\nabla f(x^\star)=0$, and \Cref{ass:setup} give $|\psi_k'(u)|\le\norm{\nabla f(z)-\nabla f(x^\star)}\le L\norm{z-x^\star}\le L(\norm{z-x^0}+R)\le L(B_{k+1}+R)$.
\end{proof}

\begin{proof}[Proof of \Cref{thm:finite}]
If $\varepsilon>LR^2$, the claim follows from the bound $f(x^0)-f^\star\le LR^2/2$ above, and no comparison is made.
We assume $0<\varepsilon\le LR^2$ in the remainder of the proof.

\pstep{Feasibility}
We first show that all operations of \Cref{alg:finite} are well defined for every sequence of reply bits, and that
\begin{equation}\label{eq:all-history-feasibility}
 x^k\in (x^0+V_k)\cap\Ball(x^0,B_k),\qquad m^k\in V_k,\qquad
 \norm{m^k}\le k\fa 0\le k\le K,
\end{equation}
where $V_0\coloneqq\{0\}$.
For $k=0$, this holds since $m^0=0$.
Suppose that~\eqref{eq:all-history-feasibility} holds for some $k<K$.
By \Cref{lem:direction}, $\widehat n^k$ is a unit vector, so $\norm{m^{k+1}}\le k+1$.
Since $B_k<B_{k+1}$, \Cref{lem:ray} gives a segment of positive length.
The line solver returns a point of this segment, and the retention step selects one of two points of the segment, so $p^{k+1}\in\Ball(x^0,B_{k+1})$.
The Gram--Schmidt procedure only divides by nonzero norms, and the span of $U$ is $V_{k+1}$, which contains $p^{k+1}-x^0$ and $m^{k+1}$.
The plane solver returns a point $y$ with $\norm y\le B_{k+1}$ (also in dimension zero, where it returns $0$), so $x^0+Uy\in(x^0+V_{k+1})\cap\Ball(x^0,B_{k+1})$.
The retention step selects either this point or $p^{k+1}$, and both lie in $(x^0+V_{k+1})\cap\Ball(x^0,B_{k+1})$.
This proves~\eqref{eq:all-history-feasibility} for $k+1$.
Note that this argument uses neither the accuracy of the directions nor the exactness of the comparisons.

\pstep{The restricted functions}
Since $U$ has orthonormal columns, the function $\varphi(y)=f(x^0+Uy)$ is convex and $L$-smooth, and, for $\norm y\le B$,
\begin{equation}\label{eq:restrictedbound}
 \norm{\nabla\varphi(y)}=\norm{U^T\nabla f(x^0+Uy)}\le\norm{\nabla f(x^0+Uy)}\le L\norm{x^0+Uy-x^\star}\le L(B+R).
\end{equation}
Together with \Cref{lem:ray}, this shows that both calls of $\SmallMin$ receive valid inputs.
Both restricted functions are defined on the whole space, so the direction probes of $\SmallMin$ may leave the ball.

\pstep{Accuracy of the subproblems and monotonicity}
As in \Cref{alg:finite}, let $C_k=(x^0+V_k)\cap\Ball(x^0,B_k)$.
The segment contains $x^k$ (at $\alpha=0$), so the retention step gives $f(p^{k+1})\le f(x^k)$, and the retained point is at least as good as the output of the line solver.
Similarly, $p^{k+1}\in C_{k+1}$, so the second retention step gives $f(x^{k+1})\le f(p^{k+1})$ and keeps the accuracy of the plane solver.
By \Cref{lem:small},
\begin{equation}\label{eq:splitaccuracy}
\begin{gathered}
 f(p^{k+1})\le\min_{0\le\alpha\le H_k}
 f(x^k-\alpha\widehat n^k)+\rho_{\rm line},\qquad
 f(x^{k+1})\le\min_{u\in C_{k+1}}f(u)+\rho_{\rm plane},\\
  f(x^{k+1})\le f(p^{k+1})\le f(x^k).
\end{gathered}
\end{equation}
In particular, the gaps $\Delta_k$ are non-increasing, and for $k\ge1$, \Cref{lem:vi} applies to $p=x^k$ on $C_k$ with $\nu=\rho_{\rm plane}$ and $B=B_k\le B_{\max}$, so that $x^k$ satisfies~\eqref{eq:vi} on $C_k$.
For $k=0$, we have $C_0=\{x^0\}$, and~\eqref{eq:vi} holds trivially.
Moreover, every iterate satisfies $\norm{x^k-x^\star}\le B_k+R\le q$.

\pstep{Accurate directions}
Let $0\le k<K$ with $\Delta_k>\varepsilon$.
By convexity, $\varepsilon<\Delta_k\le\ip{g_k}{x^k-x^\star}\le G_kq$, so $G_k>\varepsilon/q=\gamma$.
\Cref{lem:direction} then gives $\norm{\widehat n^k-n^k}\le\xi/2$, where $n^k\coloneqq g_k/G_k$.
In the next three steps, we fix such an index $k$; no estimate is needed at iterations where the gap is already below the target.

\pstep{Anchor and memory conditions}
Testing~\eqref{eq:vi} at $u=x^0$ gives $\ip{g_k}{s_k}\le\Gamma$.
Conversely, $0<\Delta_k\le\ip{g_k}{x^k-x^\star}=\ip{g_k}{s_k}+\ip{g_k}{x^0-x^\star}\le\ip{g_k}{s_k}+G_kR$.
Hence
\begin{equation}\label{eq:radial}
 -G_kR\le\ip{g_k}{s_k}\le\Gamma .
\end{equation}
If $m^k\ne0$, the points $x^0\pm B_km^k/\norm{m^k}$ belong to $C_k$, by~\eqref{eq:all-history-feasibility}.
Testing~\eqref{eq:vi} at both points and using~\eqref{eq:radial} gives $|\ip{g_k}{m^k}|\le(\Gamma-\ip{g_k}{s_k})\norm{m^k}/B_k\le(\Gamma+G_kR)\norm{m^k}/B_k$, and this also holds for $m^k=0$.
We divide by $G_k>\gamma$ and account for the error of the direction.
Since $\Gamma/\gamma=R/K\le R$, $B_k\ge16KR$, and $\xi=1/(64K)$,
\begin{equation}\label{eq:memoryapprox}
 \ip{\widehat n^k}{m^k}
 \le\frac{|\ip{g_k}{m^k}|}{G_k}+\norm{\widehat n^k-n^k}\norm{m^k}
 \le\left(\frac{R+\Gamma/\gamma}{B_k}+\xi\right)\norm{m^k}
 \le\frac9{64K}\norm{m^k}.
\end{equation}
With $\norm{m^k}\le k$ from~\eqref{eq:all-history-feasibility}, this gives the memory condition $\ip{\widehat n^k}{m^k}\le9k/(64K)=\vartheta_k$.
The anchor condition and the direction condition follow from~\eqref{eq:radial} and the previous step: $\ip{g_k}{x^k-x^0}/G_k\le\Gamma/\gamma=R/K=r_k$ and $\norm{\widehat n^k-n^k}\le\xi=\delta_k$.

\pstep{A competitor for the line search}
Consider the point $\tilde x^k\coloneqq x^k-(G_k/L)\widehat n^k$ on the ray, which is only used in the analysis.
By the triangle inequality, \eqref{eq:nonexp}, and \Cref{ass:setup}, $\norm{x^k-g_k/L-x^0}\le\norm{x^k-g_k/L-x^\star}+R\le\norm{x^k-x^\star}+R\le\norm{x^k-x^0}+2R\le B_k+2R$.
Moreover, $G_k/L\le\norm{x^k-x^\star}\le q\le20KR$, so that $\norm{\tilde x^k-(x^k-g_k/L)}=(G_k/L)\norm{\widehat n^k-n^k}\le q\xi\le20R/64$.
Hence $\norm{\tilde x^k-x^0}\le B_k+2R+20R/64<B_{k+1}$, and \Cref{lem:ray} shows that $\tilde x^k$ lies on the searched segment, that is, $0\le G_k/L\le H_k$.
By~\eqref{eq:splitaccuracy},
\begin{equation}\label{eq:splitcompetitor}
 f(x^{k+1})\le f(p^{k+1})\le f(\tilde x^k)+\rho_{\rm line}.
\end{equation}
Only the line tolerance enters here; the plane tolerance is only used for~\eqref{eq:vi}.

\pstep{Sufficient decrease}
Expanding the squared distance between the unit vectors $n^k$ and $\widehat n^k$ gives $\ip{n^k}{\widehat n^k}=1-\frac12\norm{n^k-\widehat n^k}^2\ge1-\xi^2/2$.
By~\eqref{eq:upper}, \eqref{eq:splitcompetitor}, and this bound,
\begin{equation}\label{eq:finitedescent}
 f(x^{k+1})
 \le f(x^k)-\frac{G_k^2}{L}\ip{n^k}{\widehat n^k}+\frac{G_k^2}{2L}+\rho_{\rm line}
 \le f(x^k)-\frac{1-\xi^2}{2L}G_k^2+\rho_{\rm line}
 \le f(x^k)-\frac{G_k^2}{4L},
\end{equation}
where the last step uses $\xi\le1/2$ and $\rho_{\rm line}=\gamma^2/(8L)\le G_k^2/(8L)$.

\pstep{Conclusion}
We have verified the conditions~\eqref{eq:robust-conditions} of \Cref{thm:inexact-certificate} at every index $k$ with $\Delta_k>\varepsilon$, with $u^k=\widehat n^k$ and the values~\eqref{eq:finite-certificate-budgets}.
The memory vectors are sums of unit vectors on every history, and the gaps are non-increasing by~\eqref{eq:splitaccuracy}.
Summing the error bounds gives $K+2\sum_{k<K}\vartheta_k=K+9(K-1)/64\le73K/64$ and $\sum_{k<K}(r_k+R\delta_k)=R+R/64=65R/64$.
Since $\sqrt{73/64}+65/64<3$ and $K\ge1$, this gives $S_K\le R\sqrt{73K/64}+65R/64\le3R\sqrt K$.
\Cref{thm:inexact-certificate} with $\eta=1/(4L)$ and $K\ge10\sqrt A$ therefore gives $\Delta_K\le\max\{\varepsilon,S_K^2/(\eta K^3)\}\le\max\{\varepsilon,36LR^2/K^2\}\le\max\{\varepsilon,36LR^2/(100A)\}=\varepsilon$.
Note that after the gap first drops below $\varepsilon$, the direction estimates may be inaccurate; monotonicity keeps the gap below $\varepsilon$ nevertheless.

\pstep{Number of comparisons}
The $K$ direction estimates in dimension $d$ cost $K D_d(\xi)$ comparisons, which is $O(dK\log(dK))$ by~\eqref{eq:dircost}, and the two retention steps cost $2K$ comparisons.
For the small solves, we use $A\le K^2/100$ and $q\le20KR$, which give $\gamma\ge\varepsilon/(20KR)=LR/(20AK)\ge5LR/K^3$ and $\Gamma=\gamma R/K\ge5LR^2/K^4$.
Consequently,
\begin{equation}\label{eq:nuentries}
 \rho_{\rm line}=\frac{\gamma^2}{8L}\ge\frac{25LR^2}{8K^6},\qquad
 \rho_{\rm plane}=\frac{\Gamma^2}{16LB_{\max}^2}\ge\frac{25LR^2}{5776K^{10}}\ge\frac{LR^2}{256K^{10}},
\end{equation}
where we used $16\cdot19^2=5776$ and $25/5776>1/256$.
Each solve has radius $B\le B_{\max}=19KR$ (for the line, $H_k/2\le B_{k+1}$ by \Cref{lem:ray}) and gradient bound $M=M_k\le20LKR$, so
\begin{equation}\label{eq:logratio}
 BM\le380LR^2K^2,\qquad
 \frac{BM}{\rho_{\rm line}}\le\frac{3040}{25}K^8\le122K^8,
 \qquad
 \frac{BM}{\rho_{\rm plane}}\le97280K^{12}.
\end{equation}
By \Cref{lem:small} and~\eqref{eq:logratio}, the line solve therefore costs $O(1+\log K)$ comparisons, and the plane solve costs $O(1+\log^2K)$ comparisons.
Finally, $10\sqrt A\le K\le11\sqrt A$ for $A\ge1$, which proves~\eqref{eq:cost}.
\end{proof}

\section{Unit directions from comparisons}\label{app:direction}

This appendix proves \Cref{lem:direction}.
The coordinate signs, the choice of a pivot coordinate, and the bisection of ratios in \Cref{alg:dir} follow Comparison-GDE of~\citet[Algorithm~2]{zhangli2024}.
We pad the bisection intervals by $\beta$ to absorb the error of the finite differences, so that each interval contains the true ratio.
The fixed loop lengths and the clipping of the interval endpoints keep the routine well defined, with a fixed number of comparisons, for every sequence of replies, including when the condition $\norm{\nabla F(p)}\ge\gamma$ of \Cref{lem:direction} fails.
Throughout, $F$, $n$, $L$, $\gamma$, and $\xi\in(0,1]$ are as in \Cref{lem:direction}; convexity is not used.

At a query point $p$, we write $e_i$ for the $i$-th coordinate vector and set
\begin{equation}\label{eq:dirparams}
 \tau=\frac{\gamma\xi}{32n\sqrt n},\qquad
 \underline a=\frac\gamma{2\sqrt n},\qquad
 \beta=\frac\tau{\underline a}=\frac\xi{16n},\qquad
 J=\left\lceil\log_2\frac{16\sqrt n}{\xi}\right\rceil.
\end{equation}
Here, $\tau$ bounds the error of a finite difference, $\underline a$ is a lower bound on the magnitude of the pivot coordinate when $\norm{\nabla F(p)}\ge\gamma$, and $\beta$ is the padding of the bisection intervals.
For a nonzero vector $v$, one directional comparison is
\begin{equation}\label{eq:dp}
 h_v=\frac{2\tau}{L\norm v^2},\qquad
 \mathsf{DP}(v)=\mathbf1\{F(p+h_vv)\ge F(p)\}.
\end{equation}

\begin{algorithm}[tbp]
\caption{$\Dir(p;L,\gamma,\xi)$ in dimension $n$}\label{alg:dir}
\begin{algorithmic}[1]
\Require Exact comparisons of $F$; $p,L,\gamma$, and $\xi\in(0,1]$; constants~\eqref{eq:dirparams}.
\For{$i=1,\ldots,n$} \Comment{Signs of the coordinates of $\nabla F(p)$}
\State Set $\sigma_i=+1$ if $\mathsf{DP}(e_i)=1$, and $\sigma_i=-1$ otherwise.
\EndFor
\State $j=1$.
\For{$i=2,\ldots,n$} \Comment{Running maximum: compare $i$ with current pivot $j$}
\If{$\mathsf{DP}(\sigma_i e_i-\sigma_j e_j)=1$} \State $j=i$. \EndIf
\EndFor
\State $\widehat\lambda_j=\sigma_j$.
\For{each $i\in\{1,\ldots,n\}\setminus\{j\}$} \Comment{Bisect the ratio of coordinate $i$ to the pivot}
\State $[\ell,u]=[-2,2]$.
\For{$s=1,\ldots,J$}
\State $b=(\ell+u)/2$.
\If{$\mathsf{DP}(e_i-b\sigma_j e_j)=1$}
\State $\ell=\max\{\ell,b-\beta\}$. \Comment{Padding $\beta$ absorbs the finite-difference error}
\Else \State $u=\min\{u,b+\beta\}$. \EndIf
\EndFor
\State $\widehat\lambda_i=(\ell+u)/2$.
\EndFor
\State \Return $\widehat n=\widehat\lambda/\norm{\widehat\lambda}$.
\end{algorithmic}
\end{algorithm}

\begin{proof}[Proof of \Cref{lem:direction}]
We write $g=\nabla F(p)$ and $G=\norm g$.

\pstep{Signs of directional derivatives}
For $v\ne0$, consider the finite difference $\chi_v\coloneqq(F(p+h_vv)-F(p))/h_v$, which is only used in the analysis.
By~\eqref{eq:upper}, applied to $F$, and~\eqref{eq:dp}, $|\chi_v-\ip gv|\le Lh_v\norm v^2/2=\tau$.
Since $\mathsf{DP}(v)=1$ means $\chi_v\ge0$ and $\mathsf{DP}(v)=0$ means $\chi_v<0$, the exact bit implies
\begin{equation}\label{eq:dpsign}
 \mathsf{DP}(v)=1\ \Longrightarrow\ \ip gv\ge-\tau,\qquad
 \mathsf{DP}(v)=0\ \Longrightarrow\ \ip gv\le\tau.
\end{equation}

\pstep{A positive pivot}
Let $a_i\coloneqq\sigma_ig_i$ for $1\le i\le n$ and $M_g\coloneqq\max_i|g_i|$.
By~\eqref{eq:dpsign}, a wrong sign, $a_i<0$, is only possible if $|g_i|\le\tau$.
Hence, in all cases,
\begin{equation}\label{eq:direction-sign-loss}
 a_i\ge|g_i|-2\tau\fa 1\le i\le n.
\end{equation}
Let $j_i$ denote the pivot after coordinates $1,\ldots,i$ have been processed, so that $j_1=1$.
The probe in the second loop has inner product $\ip g{\sigma_ie_i-\sigma_{j_{i-1}}e_{j_{i-1}}}=a_i-a_{j_{i-1}}$.
If the bit is one, the algorithm switches to $i$, and $a_i\ge a_{j_{i-1}}-\tau$; if the bit is zero, it keeps $j_{i-1}$, and $a_{j_{i-1}}\ge a_i-\tau$.
In both cases, $a_{j_i}\ge\max\{a_{j_{i-1}},a_i\}-\tau$ for $2\le i\le n$, and by induction, $a_{j_i}\ge\max_{t\le i}a_t-(i-1)\tau$.
For the final pivot $j=j_n$, \eqref{eq:direction-sign-loss} gives $a_j\ge\max_{1\le i\le n}a_i-(n-1)\tau\ge M_g-(n+1)\tau$.
If $G\ge\gamma$, then $M_g\ge G/\sqrt n\ge\gamma/\sqrt n$, and, since $\xi\le1$, $(n+1)\tau/M_g\le(n+1)\xi/(32n)\le1/16$, so that $a_j\ge\frac{15}{16}M_g\ge\underline a>0$.
Until the end of the normalization step, we assume $G\ge\gamma$.
Then the ratios $\lambda_i\coloneqq g_i/a_j$ satisfy $|\lambda_i|\le M_g/a_j\le16/15<2$, and $\lambda_j=\sigma_j$.
The algorithm never evaluates $a_j$.

\pstep{Bisection of ratios}
For $i\ne j$, the probe in the third loop has inner product $\ip g{e_i-b\sigma_je_j}=g_i-ba_j=a_j(\lambda_i-b)$.
By~\eqref{eq:dpsign} and $a_j\ge\underline a$, a reply of one implies $\lambda_i\ge b-\tau/a_j\ge b-\beta$, and a reply of zero implies $\lambda_i\le b+\beta$.
Hence every interval $[\ell,u]$ contains $\lambda_i$.
After a reply of one, the new width is at most $u-b+\beta$; after a reply of zero, it is at most $b+\beta-\ell$.
In both cases, the widths satisfy $W_0=4$ and $W_{s+1}\le W_s/2+\beta$ for $0\le s<J$, so that $W_J\le4\cdot2^{-J}+2\beta$.
The midpoint estimate therefore satisfies $|\widehat\lambda_i-\lambda_i|\le W_J/2\le2^{1-J}+\beta\le\xi/(8\sqrt n)+\xi/(16n)\le\xi/(4\sqrt n)$.
Since $\widehat\lambda_j=\lambda_j$, summing the squared errors over the coordinates gives $\norm{\widehat\lambda-\lambda}\le\xi/4$.

\pstep{Normalization}
For nonzero vectors $x,y$, the triangle inequality gives
\begin{equation}\label{eq:normalize}
 \left\|\frac x{\norm x}-\frac y{\norm y}\right\|
 \le\left\|\frac x{\norm x}-\frac x{\norm y}\right\|+\frac{\norm{x-y}}{\norm y}
 =\frac{|\norm y-\norm x|+\norm{x-y}}{\norm y}
 \le\frac{2\norm{x-y}}{\norm y}.
\end{equation}
Here $\norm{\widehat\lambda}\ge1$ and $\norm\lambda\ge1$, since both have the entry $\sigma_j$ in coordinate $j$.
Moreover, $\lambda/\norm\lambda=g/G$, since $a_j>0$.
Hence $\norm{\widehat n-g/G}\le2\norm{\widehat\lambda-\lambda}/\norm\lambda\le\xi/2$, which proves the accuracy claim of \Cref{lem:direction}.

\pstep{Termination for every sequence of replies}
There are $n$ sign queries, $n-1$ pivot queries, and $(n-1)J$ bisection queries, which proves~\eqref{eq:dircost}.
Every probe vector is nonzero: the pivot probes involve two different coordinates, and the bisection probes have coordinate $i$ equal to one.
For arbitrary replies, the new lower endpoint after a reply of one is at most the previous midpoint, and the new upper endpoint after a reply of zero is at least the previous midpoint, so the intervals stay nonempty.
No step divides by $a_j$, and the pivot coordinate still gives $\norm{\widehat\lambda}\ge1$.
All loops have fixed length.
For $n=1$, there are no pivot or bisection queries, and the output is $\sigma_1$.
\end{proof}

\section{Comparison solves in one and two dimensions}\label{app:small}

This appendix states \Cref{alg:small} and proves \Cref{lem:small}.
We combine \Cref{lem:direction} with the classical central-cut ellipsoid method~\citep{gls}.
Throughout, $\varphi$, $n$, $B$, $L$, $M$, $\nu$, and $C$ are as in \Cref{lem:small}.

\paragraph{The routine.}
If $M=0$ or $\nu\ge MB$, the center $0$ of $C$ is $\nu$-optimal by~\eqref{eq:lipsmall} below, and the routine returns $0$.
Otherwise, it sets
\begin{equation}\label{eq:smallparams}
 \gamma_s=\frac\nu{2B},\quad\zeta=\frac\nu{8BM},\quad
 r=\frac\nu{4M},\quad
 N=\left\lceil2n(n+1)\log\frac{B}{r}\right\rceil+1,
\end{equation}
so that $0<r<B/4$ and $0<\zeta<1/8$.
We write the current ellipsoid as $E(c,Q)=\{u:(u-c)^TQ^{-1}(u-c)\le1\}$, where $Q$ is positive definite.
For $n=2$ and a unit normal $a$, the central cut $\{u:\ip a{u-c}\le0\}$ leads to the update
\begin{equation}\label{eq:ellupdate}
 v=\frac{Qa}{\sqrt{a^TQa}},\qquad c^+=c-\frac{v}{n+1},\qquad
 Q^+=\frac{n^2}{n^2-1}\left(Q-\frac2{n+1}vv^T\right).
\end{equation}
For $n=1$, the interval update is
\begin{equation}\label{eq:intervalupdate}
 c^+=c-\frac{a\sqrt Q}{2},\qquad Q^+=\frac{Q}{4}.
\end{equation}
At a center outside $C$, the routine uses the normal $a=c/\norm c$, whose central cut keeps all of $C$.

\begin{algorithm}[tbp]
\caption{$\SmallMin(\varphi;B,L,M,\nu)$ in dimension $n\le2$}\label{alg:small}
\begin{algorithmic}[1]
\Require Convex, $L$-smooth $\varphi$ on $\R^n$, $n\in\{0,1,2\}$; exact comparisons;
$B,L>0$ and $M\ge0$ satisfying~\eqref{eq:smallbound}; accuracy $\nu>0$.
\If{$n=0$} \State \Return the unique point of the zero-dimensional space. \EndIf
\If{$M=0$ or $\nu\ge MB$} \State \Return $0$. \Comment{The center $0$ is $\nu$-optimal}\EndIf
\State Set the constants~\eqref{eq:smallparams}, and set $c=0$, $Q=B^2I_n$, $b=0$.
\For{$j=0,\ldots,N-1$} \Comment{Fixed number of cuts}
\If{$\norm c\le B$} \Comment{Feasible center}
\State Compare $\varphi(c)$ and $\varphi(b)$; retain the better point as $b$. \Comment{Best center so far}
\State $a=\Dir(c;L,\gamma_s,\zeta)$, applied to $\varphi$ in $\R^n$. \Comment{Approximate gradient direction}
\Else
\State $a=c/\norm c$. \Comment{This cut keeps all of $C$}
\EndIf
\State Update $(c,Q)$ by~\eqref{eq:ellupdate} if $n=2$, or by~\eqref{eq:intervalupdate} if $n=1$. \Comment{Keep $\{u:\ip a{u-c}\le0\}$}
\EndFor
\State \Return $b$.
\end{algorithmic}
\end{algorithm}

The direction estimates query $\varphi$ at points that may lie outside $C$, which is why \Cref{lem:small} requires $\varphi$ to be $L$-smooth on all of $\R^n$.
We first isolate the geometric part of the argument, which only uses convexity and the gradient bound~\eqref{eq:smallbound}.

\begin{lemma}[Geometry of the inner solvers]\label{lem:inner-geometry}
Let $B,M>0$, $n\in\{1,2\}$, $C=\Ball(0,B)\subset\R^n$, and let $\varphi:\R^n\to\R$ be convex and differentiable with~\eqref{eq:smallbound}.
For $0<\nu<MB$, define the target set $S=\{u\in C:\varphi(u)\le\min_{v\in C}\varphi(v)+\nu/2\}$.
Start from $c=0$ and $Q=B^2I_n$, perform $N$ updates~\eqref{eq:ellupdate} (for $n=2$) or~\eqref{eq:intervalupdate} (for $n=1$) with unit normals, using $a=c/\norm c$ at centers outside $C$, and store the best feasible center visited before each cut.
Suppose that, at every feasible center $c$,
\begin{equation}\label{eq:geometry-cut-condition}
 \varphi(c)-\min_{v\in C}\varphi(v)>\nu
 \quad\Longrightarrow\quad
 \ip a{u-c}\le0\fa u\in S.
\end{equation}
Then, with $N$ as in~\eqref{eq:smallparams}, the stored point is $\nu$-optimal on $C$.
Moreover, every update is well defined for every sequence of unit normals.
\end{lemma}

\begin{proof}
By compactness and continuity, $\varphi$ attains its minimum over $C$ at some $u^\star\in C$; we write $\varphi_C^\star=\varphi(u^\star)$.

\pstep{A ball inside the target set}
For $u,v\in C$, the segment from $v$ to $u$ lies in $C$, so the fundamental theorem of calculus and~\eqref{eq:smallbound} give
\begin{equation}\label{eq:lipsmall}
 |\varphi(u)-\varphi(v)|
 \le\int_0^1\norm{\nabla\varphi(v+t(u-v))}\norm{u-v}\,dt
 \le M\norm{u-v}\fa u,v\in C.
\end{equation}
In particular, $\varphi(0)-\varphi_C^\star\le MB$, which justifies the trivial cases of \Cref{alg:small}.
Let $r=\nu/(4M)$ as in~\eqref{eq:smallparams}, so that $r<B/4$, and let $w=(1-r/B)u^\star$.
If $\norm{u-w}\le r$, the triangle inequality gives
\begin{equation}\label{eq:small-good-distances}
 \norm u\le\left(1-\frac rB\right)\norm{u^\star}+r\le B,\qquad
 \norm{u-u^\star}\le\norm{u-w}+\frac rB\norm{u^\star}\le2r,
\end{equation}
and therefore, by~\eqref{eq:lipsmall},
\begin{equation}\label{eq:goodball}
 \varphi(u)-\varphi_C^\star\le2Mr=\frac{\nu}{2}.
\end{equation}
Hence $\Ball(w,r)\subseteq S$.

\pstep{The ellipsoid update}
Let $n=2$, suppose that the current $Q$ is positive definite, and let $T=Q^{1/2}$ be its symmetric positive definite square root.
The denominator $\sqrt{a^TQa}=\norm{Ta}$ in~\eqref{eq:ellupdate} is positive, since $a$ is a unit vector.
We set $\hat a=Ta/\norm{Ta}$ and use the coordinates $z=T^{-1}(u-c)$.
In these coordinates, the current ellipsoid is the unit ball, and since $\ip a{u-c}=\norm{Ta}\ip{\hat a}z$, the retained half is $\{z:\norm z\le1,\ \ip{\hat a}z\le0\}$.
Substituting $v=T\hat a$ into~\eqref{eq:ellupdate} gives
\begin{equation}\label{eq:ell-scaled-update}
 T^{-1}(c^+-c)=-\frac{\hat a}{n+1},\qquad
 T^{-1}Q^+T^{-1}=\frac{n^2}{n^2-1}\left(I-\frac{2}{n+1}\hat a\hat a^T\right).
\end{equation}
The new ellipsoid therefore has squared semi-axis $\frac{n^2}{n^2-1}(1-\frac{2}{n+1})=\frac{n^2}{(n+1)^2}$ in the direction $\hat a$, and squared semi-axes $\frac{n^2}{n^2-1}$ orthogonal to it.
We write a retained point as $z=s\hat a+z_\perp$ with $\ip{\hat a}{z_\perp}=0$, so that $-1\le s\le0$ and $\norm{z_\perp}^2\le1-s^2$.
Its value of the new quadratic form is at most
\begin{equation}\label{eq:ellcontain}
 \frac{(n+1)^2}{n^2}\left(s+\frac1{n+1}\right)^2
 +\frac{n^2-1}{n^2}(1-s^2)
 =1+\frac{2(n+1)}{n^2}s(s+1)\le1,
\end{equation}
since $s(s+1)\le0$ for $s\in[-1,0]$.
Hence the new ellipsoid contains the retained half.
Its squared semi-axes are positive, so $Q^+$ is positive definite.
Since the initial $Q=B^2I$ is positive definite, induction shows that every update is defined for every sequence of unit normals.
For $n=1$, we have $a\in\{-1,1\}$, and~\eqref{eq:intervalupdate} gives exactly the retained half interval, with $Q^+=Q/4>0$.

\pstep{Volume decrease}
For $n=2$, the semi-axes in~\eqref{eq:ell-scaled-update} show that each update multiplies the volume by $\sigma_2=\frac{2}{3}\cdot\frac{2}{\sqrt3}=\frac{4}{3\sqrt3}$.
Since $\mathrm e^{x}\le1/(1-x)$ for $x\in[0,1)$, we have $\mathrm e^{-1/3}\ge2/3>16/27=\sigma_2^2$, so $\sigma_2\le \mathrm e^{-1/6}$.
For $n=1$, the volume ratio is $1/2\le3/4\le \mathrm e^{-1/4}$.
In both dimensions, with $E_j$ the ellipsoid before cut $j$,
\begin{equation}\label{eq:geometry-contraction}
 \operatorname{vol}(E_{j+1})\le \mathrm e^{-1/[2(n+1)]}\operatorname{vol}(E_j)
 \fa n\in\{1,2\},\ 0\le j<N.
\end{equation}

\pstep{A $\nu$-optimal center is visited}
Suppose that every feasible center among the centers $c_0=0,c_1,\ldots,c_{N-1}$ visited by the method is worse than $\varphi_C^\star+\nu$.
At feasible centers, \eqref{eq:geometry-cut-condition} shows that the cut keeps $S$.
At an infeasible center, the normal $a=c/\norm c$ gives
\begin{equation}\label{eq:derivative-safe-geometric}
 \ip a{u-c}\le\norm u-\norm c\le B-\norm c<0\fa u\in C,
\end{equation}
so the cut keeps all of $C$.
Initially, $S\subseteq C=E_0$.
If $S\subseteq E_j$, then $S$ lies in the retained half of $E_j$, and hence in $E_{j+1}$ by the previous step.
By induction, $S\subseteq E_N$.
By~\eqref{eq:smallparams}, $N>2n(n+1)\log(B/r)$, so $\mathrm e^{-N/[2(n+1)]}<(r/B)^n$.
Together with~\eqref{eq:goodball} and~\eqref{eq:geometry-contraction}, this gives
\begin{equation}\label{eq:volumecontradiction}
 \operatorname{vol}(\Ball(w,r))
 \le\operatorname{vol}(E_N)
 \le \mathrm e^{-\frac{N}{2(n+1)}}\operatorname{vol}(\Ball(0,B))
 <\left(\frac rB\right)^n\operatorname{vol}(\Ball(0,B))=\operatorname{vol}(\Ball(w,r)),
\end{equation}
a contradiction.
Hence some visited feasible center is $\nu$-optimal, and the stored point is at least as good.
Later cuts may be inaccurate, but they cannot make the stored point worse.
\end{proof}

\begin{proof}[Proof of \Cref{lem:small}]
\pstep{Well-definedness and feasibility}
The stored point starts at $0\in C$ and is only replaced by centers in $C$, so it stays in $C$.
By \Cref{lem:direction}, the direction calls always return unit vectors, also at centers where the gradient is small, and a geometric cut is only normalized when $\norm c>B>0$.
By \Cref{lem:inner-geometry}, every update is defined for every sequence of unit normals, so \Cref{alg:small} is well defined.
By~\eqref{eq:lipsmall}, $\varphi(0)-\min_C\varphi\le MB$, so both trivial returns are $\nu$-optimal.

\pstep{Valid cuts from comparisons}
Let $c\in C$ be a center with $\varphi(c)-\varphi_C^\star>\nu$, and write $g=\nabla\varphi(c)$ and $G=\norm g$.
By convexity, $\nu<\varphi(c)-\varphi_C^\star\le\ip g{c-u^\star}\le2BG$, so $G>\nu/(2B)=\gamma_s$, and \Cref{lem:direction} (with $\zeta<1/8$) gives a unit vector $a$ with $\norm{a-g/G}\le\zeta/2$.
For $u\in S$, convexity gives $\ip g{u-c}\le\varphi(u)-\varphi(c)<-\nu/2$.
This negative margin absorbs the error in the normal.
Since $0<G\le M$ and $\norm{u-c}\le2B$,
\begin{equation}\label{eq:safecut}
 \ip a{u-c}
 =\frac1G\ip g{u-c}+\ip{a-\frac gG}{u-c}
 \le-\frac\nu{2G}+2B\zeta
 \le-\frac\nu{2M}+\frac{\nu}{4M}
 =-\frac\nu{4M}<0 .
\end{equation}
Hence the central cut keeps $S$, and~\eqref{eq:geometry-cut-condition} holds.
\Cref{lem:inner-geometry} now gives the claim of \Cref{lem:small}.

\pstep{Number of comparisons}
Each round at a feasible center uses $D_n(\zeta)$ comparisons for the direction and one comparison for the retention; rounds at infeasible centers use none.
This gives at most $N[D_n(\zeta)+1]$ comparisons.
By~\eqref{eq:smallparams} and~\eqref{eq:dircost}, both $N$ and $D_n(\zeta)$ are $O(1+\log(BM/\nu))$ for $n\le2$, since $B/r=4BM/\nu$ and $1/\zeta=8BM/\nu$.
For $n=1$, moreover, $D_1(\zeta)=1$, so the total is $2N=O(1+\log(BM/\nu))$.
In the nontrivial case, $BM/\nu>1$, so these bounds equal those stated in \Cref{lem:small}; the trivial cases use no comparisons.
Finally, for a matrix $U$ with orthonormal columns, every query of $\varphi(y)=f(x^0+Uy)$ is a single comparison of $f$.
\end{proof}

\section{Comparison cutting planes for high accuracy}\label{app:cutting}

\citet[Theorem~3]{zhangli2024} show that a cutting-plane method whose cuts are computed from comparisons finds an $\varepsilon$-optimal point with $\widetilde O(d^2)$ comparisons; see also~\citet{karabag} for an ellipsoid method on compact domains.
Since we use this bound for $A>d^2$ in \Cref{sec:algorithm}, we give a self-contained proof under \Cref{ass:setup}, for completeness.
The method runs the center-of-gravity method on the ball $\Ball(x^0,2R)$ and computes each cut from the direction estimate of \Cref{lem:direction}.
The proof uses two quantitative facts: the ball $\Ball(x^0,2R)$ contains a ball of radius $\sqrt{\varepsilon/L}$ of $\varepsilon/2$-optimal points, which bounds the number of cuts, and the gradient norm is at least $\varepsilon/(3R)$ at every center in $\Ball(x^0,2R)$ that is not $\varepsilon$-optimal, which provides the threshold $\gamma$ required by \Cref{lem:direction}.
The resulting bound $O(d^2\log A\log(dA))$ is the product of the number of cuts, $O(d\log A)$, and the number of comparisons per cut, $O(d\log(dA))$.

\begin{proposition}[Comparison cutting planes]\label{prop:cutting}
Let \Cref{ass:setup} hold, and let $0<\varepsilon\le LR^2/32$, that is, $A\ge32$.
There is a deterministic method that accesses $f$ only through the comparison oracle~\eqref{eq:oracle}, receives $x^0$, $L$, $R$, and $\varepsilon$, and returns a point satisfying~\eqref{eq:target} after $O(d^2\log A\log(dA))$ comparisons.
\end{proposition}

\begin{proof}
\pstep{A ball of good points}
We set $C=\Ball(x^0,2R)$, $r_c=\sqrt{\varepsilon/L}$, and $S=\{u\in C:f(u)\le f^\star+\varepsilon/2\}$, together with $\gamma_c=\varepsilon/(3R)$, $\xi_c=\varepsilon/(48LR^2)$, and $T_c=\lceil ed\log(2R/r_c)\rceil+1$, where $e=\exp(1)$.
Since $A\ge32$, we have $0<r_c<R$ and $0<\xi_c<1$.
If $\norm{u-x^\star}\le r_c$, then $\norm{u-x^0}\le R+r_c<2R$, and~\eqref{eq:upper} with $\nabla f(x^\star)=0$ gives $f(u)-f^\star\le L\norm{u-x^\star}^2/2\le\varepsilon/2$.
Hence $\Ball(x^\star,r_c)\subseteq S$.
The full sublevel set need not be bounded; only its part in $C$ is used.

\pstep{The method}
We use the center-of-gravity method~\citep[Section~2.1]{bubeck2015} on $C$.
Start with $P_0=C$ and the stored point $b=x^0$.
For $j=0,\ldots,T_c-1$, compute the centroid $c_j$ of $P_j$, compare $c_j$ with $b$ and retain the better point as $b$, set $a_j=\Dir(c_j;L,\gamma_c,\xi_c)$, and set $P_{j+1}=P_j\cap\{u:\ip{a_j}{u-c_j}\le0\}$.
After $T_c$ rounds, return $b$.
Computing centroids of the explicitly known sets $P_j$ is internal arithmetic, which uses no comparisons and is not counted in the oracle model.
Each $P_j$ is a convex body with nonempty interior, since a cut through the interior point $c_j$ leaves an interior on both sides; hence all centroids exist and lie in $C$.

\pstep{Cuts retain the good points}
Let $c=c_j$ be a center with $f(c)-f^\star>\varepsilon$, and write $g=\nabla f(c)$ and $G=\norm g$.
Since $\norm{c-x^\star}\le3R$, convexity and~\eqref{eq:setup} give
$\varepsilon<f(c)-f^\star\le\ip g{c-x^\star}\le3RG$ and $G=\norm{\nabla f(c)-\nabla f(x^\star)}\le3LR$.
Hence $G>\gamma_c$, and \Cref{lem:direction} gives $\norm{a_j-g/G}\le\xi_c/2$.
For $u\in S$, convexity gives $\ip g{u-c}\le f(u)-f(c)<-\varepsilon/2$, and $\norm{u-c}\le4R$, so
\begin{equation}\label{eq:cutting-safe}
 \ip{a_j}{u-c}
 \le\frac1G\ip g{u-c}+\norm{a_j-\frac gG}\norm{u-c}
 <-\frac{\varepsilon}{2G}+2R\xi_c
 \le-\frac{\varepsilon}{6LR}+\frac{\varepsilon}{24LR}<0 .
\end{equation}
Thus every cut at a center that is not $\varepsilon$-optimal retains $S$.

\pstep{Success by volume contraction}
By Gr\"unbaum's inequality~\citep[Lemma~2.2]{bubeck2015}, $\operatorname{vol}(P_{j+1})\le(1-1/e)\operatorname{vol}(P_j)\le e^{-1/e}\operatorname{vol}(P_j)$.
Suppose that none of the centers $c_0,\ldots,c_{T_c-1}$ is $\varepsilon$-optimal.
Then $S\subseteq P_{T_c}$ by the previous step, and, writing $v_d$ for the volume of the unit ball,
$v_dr_c^d\le\operatorname{vol}(S)\le\operatorname{vol}(P_{T_c})\le e^{-T_c/e}v_d(2R)^d<v_dr_c^d$, where the last step uses $T_c>ed\log(2R/r_c)$.
This contradiction shows that some center is $\varepsilon$-optimal, and the stored point is at least as good.
Later cuts may be inaccurate, but they cannot make the stored point worse.

\pstep{Number of comparisons}
Each round uses $D_d(\xi_c)$ comparisons for the direction and one for the retention, so the method uses $T_c[D_d(\xi_c)+1]$ comparisons.
Since $2R/r_c=2\sqrt A$ and $1/\xi_c=48A$, we have $T_c=O(d\log A)$ for $A\ge32$ and $D_d(\xi_c)=O(d\log(dA))$ by~\eqref{eq:dircost}, which gives $O(d^2\log A\log(dA))$.
\end{proof}

\section{Noisy comparisons}\label{app:noisy}

This section shows that the guarantee of \Cref{thm:finite} extends to comparisons that are correct with probability at least $p>1/2$ given all previous queries and replies.
Since \Cref{alg:finite} makes at most a fixed number $N$ of comparisons for every sequence of replies, including wrong ones (\Cref{lem:noisy-budget}), a union bound shows that repeating each comparison $O((p-1/2)^{-2}\log(N/\delta))$ times and taking majority votes makes all of them correct with probability at least $1-\delta$; this is the repetition argument of~\citet[Section~1.1]{feige1994} for noisy decision trees.
On this event, the run coincides with a run with exact comparisons, and \Cref{thm:finite} applies (\Cref{cor:noisy}).

\paragraph{Noise model.}
We use the model of~\citet{feige1994}, which~\citet{saha} use for comparison-based convex optimization.
When the method queries a pair $(u,v)$, it receives a bit $Y\in\{0,1\}$ such that, for a known $p\in(1/2,1]$,
\begin{equation}\label{eq:noisy-model}
 \Pr\{Y=O_f(u,v)\mid\text{all previous queries and replies}\}\ge p.
\end{equation}
This applies to every query, including repeated queries of the same pair and pairs with arbitrarily close values, and correctness refers to the tie convention in~\eqref{eq:oracle}; the probabilities may vary between queries.

\begin{lemma}[A comparison budget for every sequence of replies]\label{lem:noisy-budget}
Let $0<\varepsilon\le LR^2$, let $K$, $\xi$, $\rho_{\rm line}$, and $\rho_{\rm plane}$ be as in \Cref{alg:finite}, and set
\begin{equation}\label{eq:noisy-budget}
\begin{gathered}
 J_{\rm line}=\lceil4\log(488K^8)\rceil+1,\qquad
 J_{\rm plane}=\lceil12\log(389120K^{12})\rceil+1,\\
 J_d=\lceil\log_2(256\cdot97280K^{12})\rceil,\qquad
 N=K\bigl[D_d(\xi)+2+2J_{\rm line}+J_{\rm plane}(4+J_d)\bigr].
\end{gathered}
\end{equation}
Then, for every sequence of replies, correct or not, \Cref{alg:finite} makes at most $N$ comparisons, and $N=O(K[d\log(dK)+\log^2K])$.
\end{lemma}

\begin{proof}
\Cref{alg:finite} is well defined for every sequence of replies (\Cref{app:outer}), and each iteration makes $D_d(\xi)$ comparisons for the direction (\Cref{lem:direction}), two retention comparisons, and two calls of $\SmallMin$.
A call of $\SmallMin$ in dimension $n\in\{1,2\}$ with radius $B$, gradient bound $M$, and accuracy $\nu$ either returns without comparisons or makes at most $N_n[D_n(\zeta)+1]$ comparisons, with $N_n=\lceil2n(n+1)\log(4BM/\nu)\rceil+1$ and $\zeta=\nu/(8BM)$ by~\eqref{eq:smallparams} and the proof of \Cref{lem:small}; a call in dimension zero makes none.
For every sequence of replies, the bounds~\eqref{eq:logratio} hold, that is, $BM/\rho_{\rm line}\le122K^8$ and $BM/\rho_{\rm plane}\le97280K^{12}$.
For the ray search, $n=1$ and $D_1(\zeta)=1$, so it makes at most $2J_{\rm line}$ comparisons.
For the plane search, $N_n\le J_{\rm plane}$ for $n\in\{1,2\}$, and $D_n(\zeta)+1\le4+\lceil\log_2(128\sqrt2\,BM/\rho_{\rm plane})\rceil\le4+J_d$ by~\eqref{eq:dircost}, so it makes at most $J_{\rm plane}(4+J_d)$ comparisons.
Summing over the $K$ iterations gives the first claim, and the second follows from~\eqref{eq:dircost}.
\end{proof}

\paragraph{Majority votes.}
For $\delta\in(0,1)$, let $r$ be the smallest odd integer with $r\ge\log(N/\delta)/(2(p-1/2)^2)$; we run \Cref{alg:finite} and replace every comparison by $r$ queries of the same pair, using the majority of the $r$ replies as its result.

\begin{corollary}[Noisy comparisons]\label{cor:noisy}
Let \Cref{ass:setup} hold, let $0<\varepsilon\le LR^2$, let the replies satisfy~\eqref{eq:noisy-model}, and let $\delta\in(0,1)$.
Then \Cref{alg:finite} with majority votes returns a point satisfying~\eqref{eq:target} with probability at least $1-\delta$, and it makes at most
\begin{equation}\label{eq:noisy-cost}
 Nr\le N\left(2+\frac{2\log(N/\delta)}{(2p-1)^2}\right)
 =\widetilde O\!\left(\frac{d\sqrt A}{(p-1/2)^2}\log\frac2\delta\right)
\end{equation}
queries.
\end{corollary}

\begin{proof}
\pstep{One majority vote}
Fix the past and the next pair $(u,v)$, and let $Z_i=\mathbf1\{Y_i=O_f(u,v)\}$ for the $r$ replies $Y_1,\ldots,Y_r$.
By~\eqref{eq:noisy-model}, $Z_i=1$ with conditional probability $q_i\ge p$ given the past and $Z_1,\ldots,Z_{i-1}$.
For $U_1,\ldots,U_r$ uniform on $[0,1]$ and independent of each other and of everything else, let $B_i=Z_i\mathbf1\{U_i\le p/q_i\}$.
Then $B_i\le Z_i$, and $B_i=1$ with conditional probability $q_i\cdot p/q_i=p$ given the past, $Z_1,\ldots,Z_{i-1}$, and $U_1,\ldots,U_{i-1}$, so $B_1,\ldots,B_r$ are independent Bernoulli variables with mean $p$, independent of the past.
Since $r$ is odd, a wrong majority requires $\frac1r\sum_{i=1}^rB_i\le\frac1r\sum_{i=1}^rZ_i<\frac12$, a deviation of more than $p-1/2$ below the mean of the $B_i$.
By Hoeffding's inequality~\citep{hoeffding1963} and the choice of $r$, the conditional probability of a wrong majority is at most $\exp(-2r(p-1/2)^2)\le\delta/N$.

\pstep{All majority votes}
The method with majority votes is \Cref{alg:finite} run on the sequence of majority results, so it makes at most $N$ majority votes on every run, by \Cref{lem:noisy-budget}.
For $1\le i\le N$, let $F_i$ be the event that the $i$-th majority vote takes place and is the first wrong one.
The $i$-th pair is determined by the past before its replies are drawn, so the previous step gives $\Pr(F_i)\le\delta/N$.
The events $F_1,\ldots,F_N$ are disjoint, and some majority vote is wrong if and only if one of them occurs, so the probability that some majority vote is wrong is at most $\delta$.

\pstep{Conclusion}
If all majority votes are correct, the run coincides with the run of \Cref{alg:finite} with exact comparisons, and \Cref{thm:finite} gives~\eqref{eq:target}.
Each majority vote uses $r$ queries, and $r<2+\log(N/\delta)/(2(p-1/2)^2)$ by the choice of $r$, which gives the bound in~\eqref{eq:noisy-cost}; the last expression follows from $N=O(K[d\log(dK)+\log^2K])$ and $K=O(\sqrt A)$.
\end{proof}

\paragraph{Relation to prior work and scope.}
In the same model, \citet[Theorem~10]{saha} obtain a bound of order $dA(p-1/2)^{-2}\log(dA(p-1/2)^{-2}/\delta)$ for smooth convex functions, in expectation over the randomness of their method and without knowledge of $p$.
The model does not cover persistent errors, which repetition cannot remove~\citep{braverman2008}, or noise whose probability of a correct reply approaches $1/2$ as the two function values approach each other.

\end{document}